\documentclass[]{article}

\usepackage[]{graphicx}
\usepackage{srcltx}
\usepackage{multicol}
\usepackage{theorem}
\usepackage{rotating,array}
\usepackage{graphicx}\usepackage{color}
\usepackage{calc}
\usepackage{amsfonts,amssymb}
\usepackage{srcltx}
\usepackage{natbib}
\usepackage[utf8]{inputenc}
\usepackage{bbm}
\usepackage{mathrsfs}
\usepackage{amsmath}

\usepackage{}

\makeatletter
\def\ncite{\begingroup\def\NAT@ctype{0}\NAT@partrue\NAT@swatrue\NAT@numberstrue
\let\@cite\NAT@citenum\let\@citex\NAT@citexnum
\def\NAT@open{[}\def\NAT@close{]}\def\NAT@space{ }%
      \@ifstar{\NAT@fulltrue\NAT@cites}{\NAT@fullfalse\NAT@cites}}
\newtheorem{theorem}{Theorem}

\usepackage{nonfloat}
\usepackage[ruled,lined]{algorithm2e}

\newtheorem{lemma}{Lemma}
\newtheorem{proposition}{Proposition}	
\newtheorem{hypothesis}{Hypothesis}

\newtheorem{definition}{Definition}

\newtheorem{remark}{Remark}

\newenvironment{proof}[1][{}]
{\vspace{5pt plus 2pt minus 1pt}\par\noindent\textbf{Proof#1.} }
{\hfill$\:\:\Box$\par\vspace{8pt plus 2pt minus 1pt}}
\newenvironment{bemerkung}
{\vspace{5pt plus 2pt minus 1pt}\par\noindent\textbf{Remark. }} {\par\vspace{8pt plus 2pt
minus 1pt}} {\theorembodyfont{\upshape}}
 {\theorembodyfont{\upshape}}
 {\theorembodyfont{\upshape}}

\newcommand{\Einsbi}{{\mbox{\usefont{OT1}{cmr}{bx}{it} 1}}}

\newcommand{\N}{\mathbb N}
\newcommand{\M}{\mathbb M}

\newcommand{\R}{\mathbb R}

\newcommand{\Z}{\mathbb Z}

\newcommand{\Thetait}{\mit\Theta}
\newcommand{\Omegait}{\mit \Omega}

\newcommand\argmin{\mathop{\mathrm{argmin}}\displaylimits}

\def\ag{\left\{}
\def\ad{\right\}}

\makeatletter
\def\p@enumi{\labelenumi\expandafter\@gobble}
\def\p@enumii{\labelenumii\expandafter\@gobble}
\def\p@enumiii{\labelenumiii\expandafter\@gobble}
\def\p@enumiv{\labelenumiv\expandafter\@gobble}
\newbox\@a\newbox\@b
\newdimen\llength
\newcommand\limitsto[2]{\setbox\@a\hbox{$\SS#1$}
\setbox\@b\hbox{$\SS#2$} \ifdim\wd\@a<\wd\@b\relax\llength\wd\@b\else\llength\wd\@a\fi
\mathop{\hbox to \llength{\hskip3pt\leaders\hbox{$-\!\!\!$}\hfill
\hskip5pt\llap{$\mathord\rightarrow$}}}\limits^{#1}_{#2}} \makeatother
\renewcommand\labelenumi{{\upshape(\roman{enumi})}}

\newcommand{\fcite}[1]{\citet{#1}}
\renewcommand{\cite}{\ncite}
\def\N{\mathbb{N}}

\def\R{\mathbb{R}}

\makeatletter
\renewcommand*\descriptionlabel[1]{\hspace\labelsep
                                \normalfont\bfseries (#1)
\global\edef\@desclabel{{(#1)}}\aftergroup\@cldescdef}
\def\@cldescdef{
\global\protected@edef\@currentlabel{{\upshape\@desclabel}}}\newcounter{condition}
\newbox\cond@temp

\newenvironment{condition*}[1]{\ifnum\@listdepth>0\@warning{environment
`condition' within a list}\fi\def\labelcondition{{\rm (#1)}}\list{\labelcondition
}{\usecounter{condition}\setbox\cond@temp\hbox{#1}
\def\makelabel##1{\hss\llap{##1}}}\item}{\endlist}
\def\p@condition{\labelcondition\expandafter\@gobble}
\makeatother

\begin{document}
\title{Complexity $L^0$-Penalized $M$-Estimation: \\Consistency in More Dimensions}
\author{L. Demaret \footnote{IBB -  Institute of Biomathematics and Biometry, HMGU  Munich}
\and F. Friedrich\footnote{ETH Zentrum RZ H9, Z\"urich Switzerland; partially supported by HMGU Munich, Germany} 
\and  V. Liebscher\footnote{Ernst-Moritz-Arndt-Universit\"{a}t Greifswald, Germany}  
\and 
 G. Winkler \footnote{Ludwig-Maximilians-Universit\"at M\"unchen, Germany}
} 
\maketitle




\begin{abstract}
We study the asymptotics in $L^2$ for complexity penalized least squares regression for the discrete approximation of finite-dimensional signals on continuous domains - e.g. images - by piecewise smooth functions.

We introduce a fairly general setting which comprises most of the presently popular partitions of signal- or image- domains like interval-, wedgelet- or related partitions, as well as Delaunay triangulations. Then we prove consistency and derive convergence rates. Finally, we illustrate by way of  relevant examples that the abstract results are useful for many applications. 
\end{abstract}


\section{Introduction}\label{cons Introduction}

\begin{figure}[t]
\includegraphics[width=0.32\textwidth]{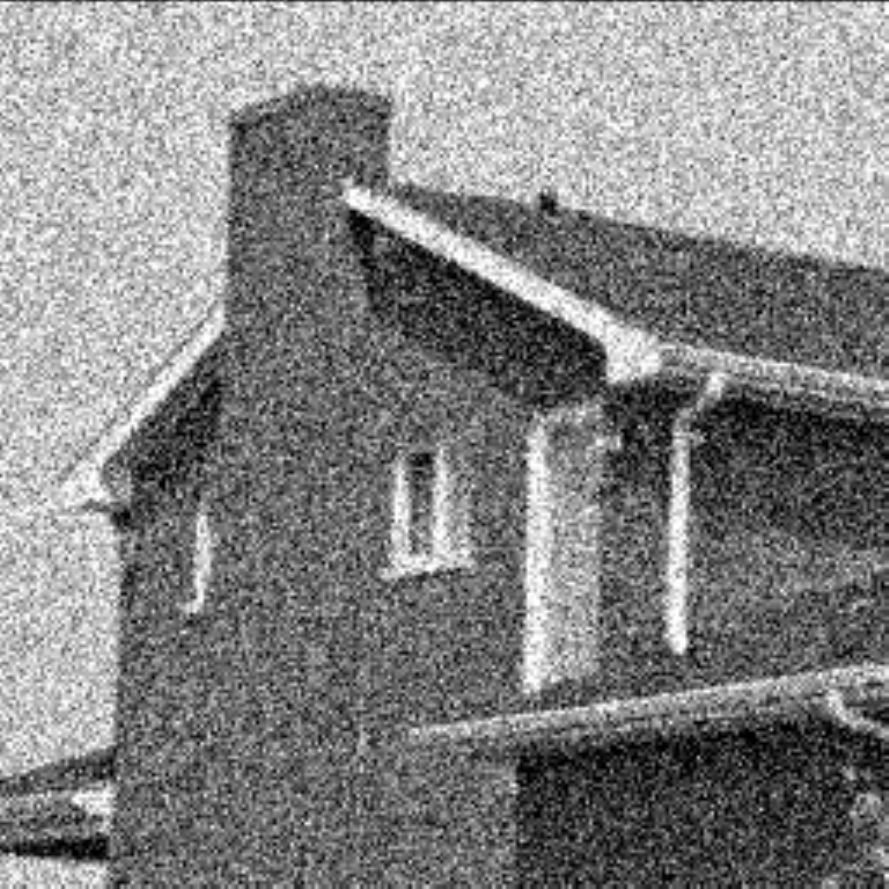}
\includegraphics[width=0.32\textwidth]{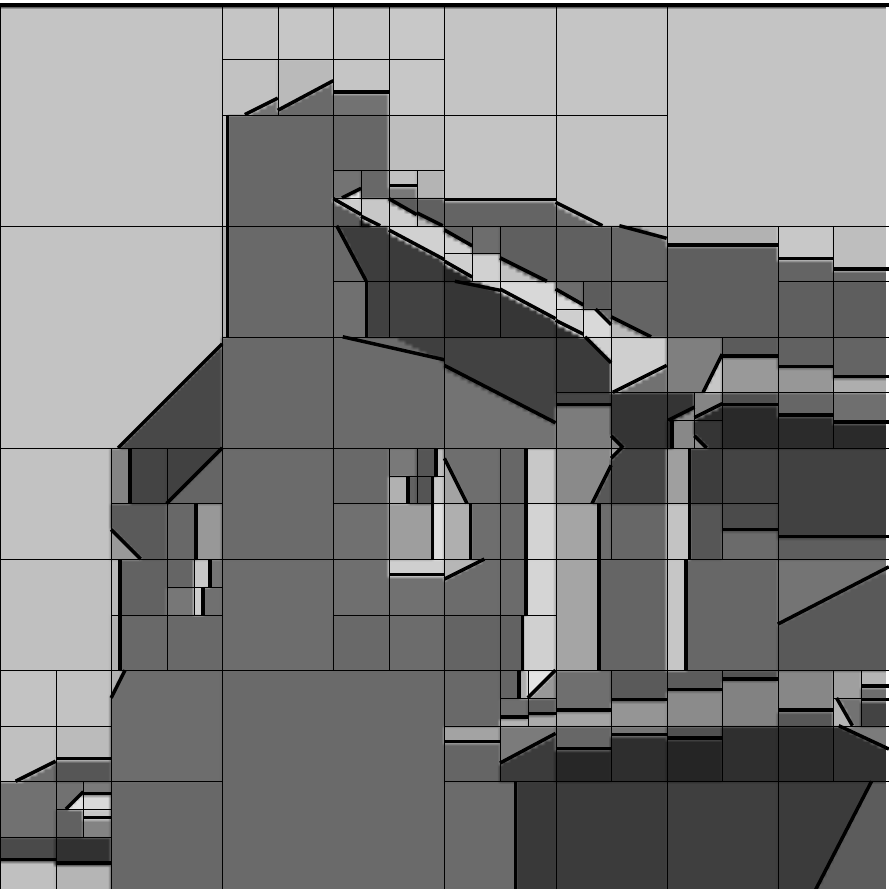}
\includegraphics[width=0.32\textwidth]{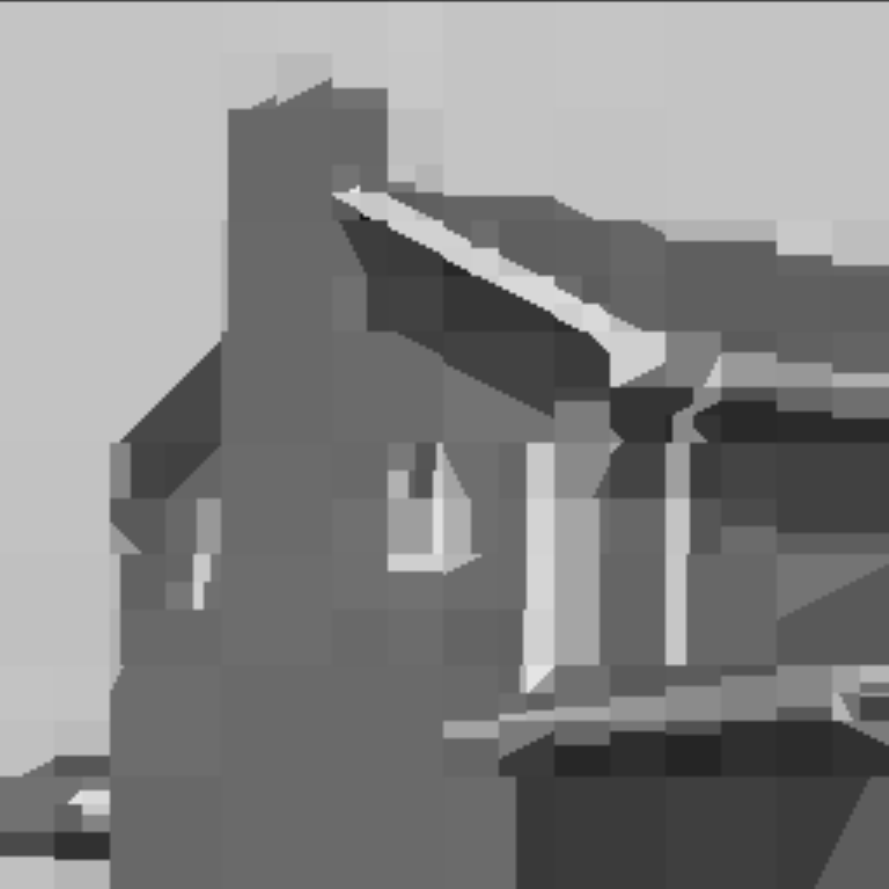}
  \caption{A noisy image (left) and (right) a fairly rough wedgelet representation for $n=256$. The (middle) picture also shows  the boundaries of the smoothness regions.}
  \label{fig:wedge}
  \end{figure}

We are going to study consistency of special complexity penalized Least Squares estimators  for noisy observations of  finite-dimensional signals on multi-dimensional domains, in particular of images. 
The estimators discussed  in the present paper are based on partitioning combined with piecewise smooth approximation. In this framework,  consistency  is proved and  convergence rates are derived in $L^2$. Finally, the abstract results are applied to a couple of relevant examples,  including popular methods like interval-, wedgelet- or related partitions, as well as Delaunay triangulations. Fig.~\ref{fig:wedge} illustrates a typical wedgelet  representation of a noisy image.

Consistency is a strong indication that an estimation procedure is meaningful. Moreover, it allows for structural insight since a sequence of discrete estimation procedures is  embedded into a common continuous setting and the quantitative behaviour of estimators can be compared. It is frequently used as a substitute or approximation for missing or vague knowledge in the real finite sample situation.   Plainly, one must be aware of various shortcomings and should not rely on asymptotics in case of small sample size. Nevertheless, consistency is a broadly accepted justification of statistical methods. Convergence rates are of particular importance, since they indicate the quality of discrete estimates or approximations and allow for comparison of different methods.

Observations or data will be governed by a simple regression model 
with additive white noise: Let $S^n=\{1,\,\ldots\,,n\}^d$ be a finite discrete signal  domain, interpreted as the discretization of the continuous domain $S^\infty=[0,1)^d$. Data $y=(y_s)_{s\in S^n}$ are available for the discrete domains at all levels $n$ and generated by the model 
\begin{equation}\label{eq regression model} 
Y_s^n=\bar{f}_s^n+\xi_s^n,\:\: n\in\mathbb{N},\:\: s\in S^n,
\end{equation}
where $(\bar{f}_s^n)_{s\in S^n}$ is a discretisation of an original or `true' signal $f$ on $S^\infty$ and $(\xi_s^n)_{s\in S^n}$ is white (sub-)Gaussian noise.

The present approach is based on a partitioning of the discrete signal domain into regions on each of which a smooth approximation of noisy data is performed. The choice of a particular  partition
is obtained by a complexity penalized least squares estimation, dependent on the data.
Between the regions, sharp breaks of intensity may happen, which  allows for edge-preserving piecewise smoothing. In one dimension, a natural way to model jumps in signals is to consider piecewise regular functions. This leads naturally to representations based on partitions consisting of intervals. The number of intervals  on a discrete line of length $n$ is of polynomial order $n^2$.

In more dimensions, however, the definition of elementary fragments is much more involved. For example, in a  discrete square of side-length $n$,  the number of all subregions is of  the exponential order $2^{n^2}$.  When dealing with images, one of the difficulties consists in constructing 
reduced sets of fragments which, at the same time, take into account the geometry of images and lead   to computationally feasible algorithms for the computation of estimators.


The estimators adopted here are minimal points of complexity penalized least squares functionals: 
 if $y=(y_s)_{s\in S^n}$ is a sample and $x=(x_s)_{s\in S^n}$ a tentative representation of $y$, the functional
\begin{equation} \label{eq:Potts1}
H^n(x,y)= \gamma|\mathscr{P}(x)|+\sum_{s\in S^n}(y_s-x_s)^2
\end{equation}
has to be minimised in $x$ given $y$; the penalty $\gamma|\mathscr{P}(x)|$ is the number of subdomains into which the entire domain is divided and on which $x$ is smooth in a sense to be made precise by the choice of suitable function spaces (see Sections \ref{suse The Regression Model} and \ref{se Special Segmentations}); $\gamma$ is a tuning parameter.
  This automatically results in a sparse representation  of the function. Due to the non-convexity of this penalty one has to solve hard optimisation problems.

These are not computationally feasible, if all possible partitions of the signal domain are admitted. 
A most popular attempt to  circumvent this nuisance is simulated annealing, see for instance the seminal paper 
\citet{GemanGeman:84}. This paper  had a considerable  impact on imaging; the authors transferred models from statistical physics to image analysis as prior distributions in the framework of Bayesian statistics. This approach was intimately connected with Markov Chain Monte Carlo Methods  like Metropolis Sampling and Simulated Annealing, cf. \citet{Winkler:2003}.

On the other hand, transferring spatial complexity to time complexity like in such metaheuristics, does not remove the basic problem; it rather transforms it.  Such algorithms  are not guaranteed to find the optimum or even a satisfactory near-optimal solution, cf. \citet{Winkler:2003}, Section 6.2. All metaheuristics will eventually encounter problems on which they perform poorly. Moreover, if the number of partitions grows at least exponentially,  it is difficult to derive useful uniform bounds on the projections of noise onto the subspaces induced by the partitions.
Reducing the search space drastically allows to design exact and fast algorithms.  Such a  reduction basically amounts  to  restrictions on   admissible partitions of the signal domain. There are various suggestions, some of them mentioned initially.

In one dimension, regression onto piecewise constant functions was proposed by the legendary \citet{Tukey:1961} who called respective representations regressograms. The functional \eqref{eq:Potts1} is by some (including the authors) referred to as the \textit{Potts functional}. It was introduced in \citet{Potts:52} as a generalization of the well-known Ising model, \citet{Ising:25}, from statistical physics from two to more spins. It was suggested by \citet{Lenz:1920} and  penalizes the length of contours between regions of constant spins. In fact, in one dimension a partition $\mathscr{P}$ into say $k$ intervals on which the signal is constant admits $k-1$ jumps and therefore has contour-length $k-1$.

The one-dimensional Potts model for signals was studied in detail in a series of theses and articles, see
\citet{WiLi:2002,WinklerKempeLiebscherWittich:2003,liebwi:99,Kempe:2003,Friedrich:2005,
WinklerWittichKempeLiebscher:2005,WittichKempeWinklerLiebscher:2008,FriedrichKempeLiebscherWinklerLiebscher:2007}. 
Consistency was first adressed in \citet{Kempe:2003} and later on exhaustively treated in \citet{LiebscherBoysenMunk:2009} and \citet{BoysenLiebscherMunkWittich:2007}. Partitions consist there of intervals. Our study of the multi-dimensional case started with the thesis \citet{Friedrich:2005}, see also \citet{Friedrich:2007}.

In two  or more dimensions, the model \eqref{eq:Potts1} differs substantially from the classical Potts model. The latter penalizes the \textit{length of contours} - locations of intensity breaks - whereas \eqref{eq:Potts1} penalizes the \textit{number of  regions}. This allows for instance to perform well on  filamentous structures,  albeit they have long borders compared to their area. 

Let us give an informal introduction into the setting. The aim is   to estimate a function $f$ on the $d$-dimensional unit cube $S^\infty=[0,1)^d$ from discrete data. 
To this end, $S^\infty$ and $f$ are discretized to cubic grids $S^n=\{1,\,\ldots\,,n\}^d$, $n\in\mathbb{N}$, and  functions $\bar{f}^n$ on $S^n$. On each stage $n$, data $y_s^n$, $s\in S^n$, is available, i.e. noisy observations of the $\bar{f}_s^n$. We will prove  $L^2$-convergence of complexity penalized  least squares estimators  $\hat{f}^n(y)$ 
(Section \ref{Complexity Penalized M-Estimation}) for the $\bar{f}^n$ (Section \ref{suse The Regression Model}) to $f$ and derive convergence rates, first in the general setting. We are faced with three kinds of error: the error caused by noise, the approximation and the often ignored  error. Noise is essentially controlled regardless of the specific form of $f$. For the approximation and the discretisation error special assumptions on the function classes in question are needed.

Because of the approximation error term, there are deep connections to approximation theory. In particular, when dealing with piecewise regular images,  non linear approximation rates obtained by wavelet shrinkage methods are known to be suboptimal, as discussed in
\citet{KorostelevTsybakov:1993} or \citet{Donoho:1999}.	 In the last decade, the challenging problem to improve upon wavelets has been addressed in very different directions.  

The search for a good paradigm for detecting and representing curvilinear discontinuities of 
bivariate functions remains a fundamental issue in image analysis. Ideally, an efficient representation 
should use atomic decompositions which are local in space (like wavelets), but also possess appropriate directional 
properties (unlike wavelets). One of the most prominent examples is given by
curvelet representations, which are based on multiscale directional filtering combined with 
anisotropic scaling. \citet{CandesDonoho:2002}  proved that thresholding of curvelet coefficients 
provides estimators which yield the minimax convergence rate up to a logarithmic factor for piecewise $\mathscr{C}^2$ functions with $\mathscr{C}^2$ boundaries.  Another interesting representation is given by bandelets as proposed in \citet{LePennecMallat:2005}.
Bandelets are based on optimal local warping in the image domain relatively to the geometrical flow and 
  \citet{DossalMallatLePennec:2011} proved  also optimality of the minimax convergence rates of their bandelet-based estimator, for a larger class of functions  including  piecewise  $\mathscr{C}^\alpha$ functions with $\mathscr{C}^\alpha$ boundaries.

The bidimensional examples discussed in Section \ref{se Special Segmentations} are based on more  geometrical constructions, to which the abstract framework proposed in Section 
\ref{se Consistency} applies.
  
Wedgelet partitions were  introduced  by \citet{Donoho:1999} and 
belong to the class of shape-preserving image segmentation methods. 
The decompositions are based on local polynomial approximation on some adaptively 
selected leaves of a quadtree structure. 
 The use of a suitable data structure allowed for the development of fast algorithms for 
wedgelet decomposition, see \citet{Friedrich:2007}.  

An alternative is provided by anisotropic Delaunay triangulations, which have been proposed in the context of image compression in \citet{DemaretDynIske:2006}.
The flexible design of the representing system allows for a particularly fine selection of triangles fitting   the   anisotropic geometrical features of images. In contrast to curvelets, such representations preserve the advantage of wavelets and are still able to approximate point singularities optimally, see \citet{DemaretIske:2011}.

Both wedgelet representations and anisotropic Delaunay triangulations lead to optimal non linear approximation rates for some classes of piecewise smooth functions.  In the present paper, we prove
 optimality  also for the convergence rates of the estimators.   More precisely, we prove  strong consistency rates of                                                                                                                    
 $$
   O(\varepsilon_n^	{2 \alpha/(\alpha+1)} \log (\varepsilon_n)), \varepsilon_n=\sigma^2/n^d,
 $$  
where $\sigma^2$ is the variance of noise and $\alpha$ is a parameter controlling piecewise regularity. Such rates  are known to be optimal up to the logarithmic factor. 

\citet{birge2007minimal} showed recently  that, in a similar setting, optimal rates without the $\log$ factor may be achieved with penalties different from those in  \eqref{eq:Potts1},  and not merely proportional to the number of pieces.  In the present work, we explicitly restrict our attention to the classical penalty given by the number of pieces as in \eqref{eq:Potts1}, noting that this corresponds to an \emph{ansatz} which is currently  popular in the signal community. One of the main reasons  is  the connection to sparsity. The generalization  of the proofs in this paper is straightforward but would be rather technical an thus might obscure the main ideas.

We address first noise and its projections to the approximation spaces, see Section \ref{se Noise and its Projections}.  In Section \ref{se Consistency},  we derive convergence rates in the general context. Finally, in Section \ref{se Special Segmentations}, we illustrate the abstract results by specific applications. 
 Dimension $1$ is included, thus generalising the results from \citet{LiebscherBoysenMunk:2009} to piecewise polynomial regression and piecewise Sobolev classes.
Our two-dimensional examples, wedgelets and Delaunay triangulations, both rely on a geometric and edge-preserving representation. Our main motivation are the optimal approximation properties of these methods,  the key feature to apply the previous framework being an appropriate discretization of these schemes.  


\section{The Setting} In this section we introduce the formal framework for piecewise smooth representations, the regression model for data, and the estimation procedure.
\subsection{Regression and Segmentations}\label{suse The Regression Model}

Image domains will generically be denoted by $S$. We choose   $S^\infty=[0,1)^d$, $d\in\N$, as the continuous and $S^n=\{1,\,\ldots\,,n\}^d$ as the generic discrete image domain. Let $f\in L^2(S^\infty)$ represent the `true' image which has to be reconstructed from noisy discrete data. For the latter, we adopt a simple linear regression model of the form
\begin{equation}\label{eq regression model} 
Y_s^n=\bar{f}_s^n+\xi_s^n,\:\: n\in\mathbb{N},\:\: s\in S^n.
\end{equation}
The noise variables $\xi_s^n$ in the regression model 
are random variables on a common probability space $(\Omegait,\mathscr{F},\mathbb{P})$. $\bar{f^n}=(\bar{f^n}_s)_{s\in S^n}$ is a discretisation of $f$. To be definite,
divide $S^\infty$ into $n^d$ semi-open cubes
\begin{displaymath}
I_{i_1,\ldots,i_d}^{n}=\prod_{1\le j\le d}[(i_j-1)/n,i_j/n),\:1\le i_j\le n,
\end{displaymath}
of volume $1/n^d$ and for $g\in L^2(S^\infty)$ take local means
\begin{displaymath}
\bar{g}_s^n=n^d\int_{I_s} g(u)\,du,\:\:s\in S^n.
\end{displaymath}
This specifies maps $\delta^n$ from $L^2(S^\infty)$ to $\mathbb{R}^{S^n}$ by
\begin{equation} \label{eq:delta}
\delta^n g= (\bar{g}_s^n)_{s\in S^n}.
\end{equation}
Conversely, embeddings of $\mathbb{R}^{S^n}$ into $L^2(S^\infty)$ are defined by
\begin{equation} \label{eq:iota} z=(z_s)_{s\in S^n}\longmapsto 
\iota^n z=\sum_{s\in S^n} z_s\Einsbi_{I_s}.	
\end{equation}
As an aid to memory, keep the following  chain of maps in mind:
\begin{displaymath}
L^2(S^\infty)\stackrel{\delta^n}{\longrightarrow}\mathbb{R}^{S^n}\stackrel{\iota^n}{\longrightarrow}L^2(S^\infty).
\end{displaymath}


%

In absence of noise, $f$ is approximated by the functions $\iota^n \bar{f}^n=\iota^n \delta^n f$ in any precision. The main task thus will be  to control noise.  In fact, the function $\iota^n\delta^n f=\iota^n \bar{f}^n$ is the conditional expectation of $f$ w.r.t. the ($\sigma$)-algebra $\mathscr{A}^n$ generated by the cubes $I_s^n$ and convergence can be seen by a martingale argument.

We are dealing with estimates of $f$ or rather of $\bar{f}^n$ on each level $n$. 
An  image domain $S$ will be partitioned by the method into sets, on which the future representations are members of initially chosen  spaces of smooth functions. To keep control, we choose  a class $\mathscr{R}\subset 2^S$ of \textit{admissible fragments} and later on, these will be rectangles, wedges or triangles. A subset $\mathscr{P}\subset 2^S$ is a \textit{partition} if  (a) the elements in $\mathscr{P}$ are mutually disjoint, and (b) $S$ is the union of all $P\in\mathscr{P}$. We will only consider partitions  $\mathscr{P}\subset\mathscr{R}$.
In addition, we choose a subset $\mathfrak{P}$ of all  partitions and call its elements  {\em admissible partitions}.

For each fragment $P\in\mathscr{R}$, we choose a finite dimensional  linear space $\mathscr{F}_P$ of real functions on $S$ which vanish off $P$.  Examples are spaces of constant functions or polynomials of higher degree. This space is determined by the maximal local smoothness of $f$. If  $\mathscr{P}\in\mathfrak{P}$ and $f_\mathscr{P}=(f_P)_{P\in\mathscr{P}}$ is  a family of such functions, we also denote by $f_\mathscr{P}$ the function defined on all of $S$ and whose restriction to $P$ is equal to $f_P$ for each $P \in \mathscr{P}$.
The pair $(\mathscr{P},f_\mathscr{P})$ is a \textit{segmentation}
and each element $(P,f_P)$ is a  \textit{segment}. 

For each partition $\mathscr{P}$ , define the linear space $\mathscr{F}_{\mathscr{P}} =
\mathrm{span} \{\mathscr{F}_P : P \in \mathscr{P}\}$.  A family 
of segmentations is called a \emph{segmentation class}. In particular, let
$$
   \mathfrak{S}(\mathfrak{P},\mathfrak{F}) := \left\{ (\mathscr{P},f) : \mathcal{P} \in \mathfrak{P}, f \in \mathscr{F}_{\mathscr{P}}\right\}
$$
with partitions in $\mathfrak{P}$ and functions in $\mathfrak{F}= \{ \mathscr{F}_{\mathscr{P}} : 
\mathscr{P}\in\mathfrak{P}\}$. 
\subsection{Complexity Penalized Least Squares Estimation}\label{Complexity Penalized M-Estimation}
We want to produce appropriate discrete representations or estimates of the underlying function $f$ on the basis of random data $Y$ from the regression model (\ref{eq regression model}). We are watching out for  a segmentation  which is in proper balance between fidelity to data and complexity.

We decide in advance on a class $\mathfrak{S}$ of (admissible) segmentations which should contain the desired representations.  The segmentations, given data $Y^n$, are scored by the functional 	
\begin{equation}\label{eq Potts functional}
H_\gamma^n: \mathfrak{S}^n\times\mathbb{R}^{S^n}\longrightarrow\mathbb{R},\: H_\gamma^n \left((\mathscr{P},f_\mathscr{P}),Y^n\right) = \gamma|\mathscr{P}| +
\|f_\mathscr{P}-Y^n\|^2,
\end{equation}
with $\gamma\ge 0$ and $|\mathscr{P}|$ the cardinality  of $\mathscr{P}$. The symbol $\|\,\cdot\,\|$ denotes the $\ell^2$-norm on $\R^{S^n}$. The last term measures fidelity to data. 
The other term is a rough measure of overall  smoothness.  As estimators for $f$ given data $Y$ we choose minimisers 
$(\hat{\mathscr{P}^n},\hat{f}^n)$ of (\ref{eq Potts functional}). Note that both $\hat{\mathscr{P}^n}$ and $\hat{f}^n$ are random since $Y^n$ is random.

The definition makes sense since minimal points of (\ref{eq Potts functional}) do always exist. This can easily verified by the \emph{reduction principle}, which relies on the decomposition
\begin{displaymath}
\min_{\mathscr{P}\in\mathfrak{P}^n,f_{\mathscr{P}}\in\mathscr{F}_\mathscr{P}} H_\gamma^n((\mathscr{P},f_{\mathscr{P}}), Y)=
\min_{\mathscr{P}\in\mathfrak{P}^n}\left(\gamma|\mathscr{P}|+\min_{f_\mathscr{P}\in\mathscr{F}_\mathscr{P}}\|f_\mathscr{P}-Y\|^2\right).
\end{displaymath}
Given $\mathscr{P}$, the inner minimisation problem has as unique solution the orthogonal projection $\hat{f}^n_\mathscr{P}$ of $Y$ to $\mathscr{F}_\mathscr{P}$. The outer minimisation problem is finite and hence a minimum of (\ref{eq Potts functional}) exists.  Let us pick one of the minimal points  $\hat{f}^n$.

\section{Noise and its Projections}\label{se Noise and its Projections}
For consistency, resolutions at infinitely many  levels  are considered simultaneously. Frequently, segmentations are not defined for all $n\in\N$ but only for a cofinal subset of $\N$. Typical examples are all dyadic partitions like quad-trees or dyadic wedgelet segmentations where only indices of the form $n=2^p$ appear. Therefore we adopt the following convention:
\begin{center}
\hspace{0ex}The symbol $\M$ denotes any infinite subset of $\N$
endowed with the natural order $\le$. 
\end{center}
$(\M,\le)$ is a totally ordered set and we may consider nets $(x_n)_{n\in\M}$. For example $x_n\rightarrow x$, $n\in\M$, means that $x_n$ convergences to $x$ along $\M$.  We deal similarly with notions like $\limsup$ etc. Plainly, we might resort to subsequences instead but this would cause a change of indices which is notationally inconvenient. 

\subsection{Sub-Gaussian Noise and a Tail Estimate}\label{se Subgaussian Noise}
We introduce now the main hypotheses on noise accompanied by a brief discussion. 
The core of the arguments in later sections is the  tail estimate (\ref{eq MainTailEstimate}) below.

As Theorem \ref{hypsubgaussian} will show, the appropriate framework are   
\emph{sub-Gaussian} random variables. A random variable $\xi$ enjoys this property if one of the following conditions is fulfilled:

\begin{theorem}
\label{th_subgaussian_equivalence}
The following two conditions on a random variable $\xi$ are equivalent:

\noindent (a) There is $a\in\mathbb{R}$ such that 
\begin{equation}\label{eq subgauss def}
\mathbb{E}(\exp{(t\xi)})\le \exp(a^2 t^2/2) \mbox{ for } t >0
\end{equation}
(b) $\xi$ is centred and \emph{majorized in distribution} by some centred Gaussian variable $\eta$, i.e.
\begin{displaymath}
\hbox{there is  }c_0\ge 0\hbox{  such that  }\mathbb{P}(|\xi|\ge c)\le\mathbb{P}(|\eta|\ge c)\hbox{  for all  }c> c_0.
\end{displaymath}
\end{theorem}
This and most other facts about sub-Gaussian variables quoted in this paper are verified in the first few sections of the monograph \fcite{BuldyginKozachenko:2000}; one may also consult \fcite{Petrov:1975}, Section III.\S 4.

The definition in \emph{(a)} was given in the celebrated paper \fcite{Chow:1966} which uses the term \textit{generalized Gaussian variables}.
 The closely related concept of semi-Gaussian variables - which requires symmetry of $\xi$ - seems to go back to \fcite{Kahane:1963}. 

The class of all sub-Gaussian random variables living on a common probability space $(\Omegait,\mathscr{A},\mathbb{P})$ is denoted by 
$\hbox{Sub}(\Omegait)$. The \emph{sub-Gaussian standard} is the number
\begin{displaymath}\label{eq DefTau}
\tau(\eta)=\inf\{a\ge0: a\:\:\hbox{is feasible in}\:\:(\ref{eq subgauss def})\}.
\end{displaymath}
The infimum is attained and hence is a minimum.
$\hbox{Sub}(\Omegait)$ is a linear space, $\tau$ is a norm on $\hbox{Sub}(\Omegait)$  if variables differing on a null-set only are identified. 
$(\hbox{Sub}(\Omegait),\tau)$ is a Banach space. It is important to note that $\hbox{Sub}(\Omegait)$ is strictly contained in all spaces $L^p_0(\Omegait)$, $p\ge 1$, the spaces of all centred variables with finite $p^{th}$ order absolute moments. 
\begin{remark} 
The most prominent  sub-Gaussians are centred Gaussian variables $\eta$ with  standard deviation $\sigma$ and $\tau(\eta)=\sigma$. For them inequality (\ref{eq subgauss def}) is an equality with $a=\sigma$. 
The specific characteristic of sub-Gaussian variables are tails lighter than those of Gaussians, as expressed in \emph{(b)} of Theorem \ref{th_subgaussian_equivalence}. 
\end{remark}

The following theorem is essential in the present context.

\begin{theorem}\label{hypsubgaussian}                    
For each $n\in \mathbb{M}$, suppose that the variables $\xi_s^n$, $s\in S^n$, are independent. Then 

\noindent (a) Suppose that there is a real number $\beta >0$ such that 
for each $n\in\mathbb{M}$ and real numbers $\mu_s$, $s\in S^n$, and each $c\in\mathbb{R}_+$, the estimate
\begin{equation}\label{eq MainTailEstimate}
\mathbb{P}\left(\left|\sum_{s\in S_n}\mu_s\xi_s^n\right|\ge c\right)\le 2\cdot \exp\left(-\displaystyle \frac{c^2}{\beta\sum_{s \in S_n}\mu_s^2}\right)
\end{equation}
holds.
Then all variables $\xi_s^n$ are sub-Gaussian with a common scale factor $\beta$.

\noindent (b) Let all variables $\xi_s^n$ be sub-Gaussian. Suppose further that  
\begin{equation}\label{eq defBeta}
\beta=2\cdot\sup \{\tau^2(\xi_s^n): n\in\mathbb{M},\:s\in S^n\}<\infty.
\end{equation}
Then (a) is fulfilled with  this factor $\beta$.
\end{theorem}
This is probably folklore. On the other hand, the proof is not straightforward and therefore we supply it in an Appendix.

\begin{remark} \label{white gaussian noise}\rm For  white Gaussian noise 
one has  $\tau(\xi_s^n) = \sigma$ and hence $\beta = 2 \sigma^2$. 
\end{remark}

\subsection{Noise Projections}

In this section, we quantify projections of noise. Choose for each $n \in \mathbb{M}$ a class $\mathscr{R}^n \subset 2^{S^n}$ of admissible segments  over $S^n$ and a set $\mathfrak{P}^n$ of admissible partitions.  As previously,  for each $P \in \mathscr{R}^n$,
a linear function space $\mathscr{F}_P$ is given. 
We shall denote orthogonal $L^2$-projections onto the  linear spaces $\mathscr{F}_\mathscr{P}=\mathrm{span}\{\mathscr{F}_P:P\in\mathscr{P}\}$
by $\pi_{\mathscr{P}}$. 

The following result provides $L^2$-estimates for the projections of noise to these spaces, as there are  more and more admissible segments.

 \begin{proposition}\label{prop estimate projections}
 Suppose that $\mathrm{dim} \,\mathscr{F}_P\le D$ for all $n\in\mathbb{M}$ and each $P\in\mathscr{R}^n$.
 Assume in addition that there is a number $M>0$ such that for some $\kappa>0$
\begin{displaymath}
\left|\mathscr{R}^n\right|\ge M\cdot n^\kappa\:\:\hbox{eventually}.
\end{displaymath}
Then for each $C>(1/\kappa+1)\beta D$ and for almost all $\omega\in\Omegait$
\begin{displaymath}
\left\|\pi_{\mathscr{P}^n}\xi^n(\omega)\right\|^2\le C|\mathscr{P}^n|\ln(|\mathscr{R}^n|) \: \;\:\hbox{ for eventually all}\:\:n\in\mathbb{M}\:\:
\hbox{and each}\:\:\mathscr{P}^n\in\mathfrak{P}^n.
\end{displaymath}
 \end{proposition}
%

This will be proven at a more  abstract level. No structure of the finite sets $S^n$ is required. Nevertheless, we adopt all definitions from Section \ref{cons Introduction} \textit{mutatis mutandis}.  All Euclidean spaces $\mathbb{R}^k$ will be endowed with their natural inner products $\langle\,\cdot,\,\cdot\,\rangle$ and respective norms. Projections onto linear subspaces $\mathscr{H}$ will be denoted by $\pi_\mathscr{H}$.

\begin{theorem}\label{th Maximal equality}
Suppose that the noise variables $\xi_s^n$ fulfill (\ref{eq MainTailEstimate}) accordingly. 
Consider finite nonempty collections $\mathfrak{H}^n$ of linear subspaces in $\mathbb{R}^{S^n}$ and assume that the dimensions of all  subspaces $\mathscr{H}\in\mathfrak{H}^n$, $n\in\mathbb{M}$, are uniformly  bounded by some number $D\in\mathbb{N}$. 
Assume in addition that there is a number $M>0$ such that for some $\kappa>0$
\begin{equation*}
\left|\mathfrak{H}^n\right|\ge M\cdot n^\kappa\:\:\hbox{eventually}.
\end{equation*}
Then for each $C>(1/\kappa+1)\beta D$ and for almost all $\omega\in\Omegait$
\begin{equation*}
\left\|\pi_\mathscr{H}\xi^n(\omega)\right\|^2\le C\ln(|\mathfrak{H}^n|) \quad \hbox{for eventually all}\:\:n\in\mathbb{M},\;\hbox{and each }\:\:\mathscr{H}\in\mathfrak{H}^n\;.
\end{equation*}
\end{theorem}

Note that $\|\cdot\|$ is Euclidean norm in the spaces $\mathbb{R}^{S_n}$, since each $\xi^n(\omega)$ is simply a vector. The assumption in the theorem can be reformulated as $|\mathfrak{H}^n|^{-1}=O(n^{-\kappa})|$.

\begin{proof}
Choose $n\in\mathbb{M}$ and $\mathscr{H}\in\mathfrak{H}^n$ with $\dim\mathscr{H}=d_n$. Let $e_i$, $1\le i\le d_n$ be some orthonormal basis of $\mathscr{H}$. Observe that for any real number $c>0$, 
\begin{equation*}
\sum_{i=1}^{d_n}|\langle\xi^n(\omega),e_i\rangle|^2>c^2\ln|\mathfrak{H}^n|\:
\end{equation*}
implies that
\begin{equation*}|\langle\xi^n(\omega),e_i\rangle|^2>\frac{c^2}{d_n}\ln|\mathfrak{H}^n|\:\:\hbox{for at least one }\:i=1,\,\ldots\,,d_n.
\end{equation*}
We derive a series of inequalities:
\begin{eqnarray*}
&&\mathbb{P}\left(\left\| \pi_\mathscr{H} \xi^n\right\|^2>c^2\ln |\mathfrak{H}^n|\right)=
\mathbb{P}\left(\sum_{i=1}^{d_n}|\langle \xi^n,e_i\rangle|^2 >c^2\ln|\mathfrak{H}^n|\right)\\
&\le&\mathbb{P}\left(  
\bigcup_{i=1}^{d_n} \{|\langle \xi^n,e_i\rangle|^2>\frac{c^2}{d_n}\ln{|\mathfrak{H}^n}|\}
\right)
\le\sum_{i=1}^{d_n}\mathbb{P}\left(|\langle\xi^n,e_i\rangle|^2>\frac{c^2}{d_n}\ln{|\mathfrak{H}^n|}\right)\\
&=&\sum_{i=1}^{d_n}\mathbb{P}\left(\left|\sum_{s\in S^n}\xi^n_s e_{i,s}\right|>c\left(\ln{|\mathfrak{H}^n|}/d_n\right)^{1/2}\right),\\
\end{eqnarray*}
where the first inequality holds because of the introductory implication. By (\ref{eq MainTailEstimate}) we may continue with
\begin{displaymath}
\le2\cdot d_n
\exp\left(\frac{-c^2\ln{|\mathfrak{H}^n|}}{\beta d_n \sum_{s\in S^n}e_{i,s}^2}\right)
\le2\cdot D \cdot|\mathfrak{H}^n|^{\frac{-c^2}{\beta D} }.
\end{displaymath}
Therefore
\begin{eqnarray*}
&&\sum_{n\in\mathbb{M}, \mathscr{H}\in\mathfrak{H}^n}\mathbb{P}\left(\left\| \pi_\mathscr{H} \xi^n\right\|^2>c^2\ln |\mathfrak{H}^n|\right)
\le 2D\sum_{n\in\mathbb{M}, \mathscr{H}\in\mathfrak{H}^n} |\mathfrak{H}^n|^\frac{-c^2}{\beta D}
\le 2D\sum_{n\in\mathbb{M}}|\mathfrak{H}^n||\mathfrak{H}^n|^\frac{-c^2}{\beta D}\\
&\le& 2D\sum_{n\in\mathbb{M}} \left(\frac{1}{M}\cdot n^{-\kappa}\right)^{\frac{c^2}{\beta D}-1}=2D\cdot M^{1-c^2/(\beta D)}\sum_{n\in\mathbb{M}}n^{-\kappa(\frac{c^2}{\beta D}-1)}.
\end{eqnarray*}
For  $C=c^2>(1/\kappa+1)\beta D$ the negative exponent becomes larger than 1 and the sum becomes finite. Enumeration of each $\mathfrak{H}^n$ and subsequent concatenation yields a sequence of events. The Borel-Cantelli lemma yields
\begin{displaymath}
\mathbb{P}(\|\pi_\mathscr{H}\xi^n\|> C\ln |\mathfrak{H}_n|\quad\hbox{ for finitely many } (n,\mathscr{H}) 
\hbox{ with } \mathscr{H} \in \mathfrak{H}^n)=1.
\end{displaymath}
This implies the assertion.
\end{proof}

Now let us now prove   the desired result.
\begin{proof}[ Proof of Proposition \ref{prop estimate projections}]
We apply Theorem \ref{th Maximal equality} to the collections $\mathfrak{H}^n=\{\mathscr{F}^n_R:R\in\mathscr{R}^n\}$. Then $|\mathfrak{H}^n|=|\mathscr{R}^n|$. Since for each  $\mathscr{P}^n \in \mathfrak{P}^n$ the spaces $\mathscr{F}^n_P$, $P\in\mathscr{P}^n$, are mutually orthogonal, one has  for $z\in\mathbb{R}^{S^n}$ that 
\begin{displaymath}
\|\pi_{\mathscr{P}^n} z\|^2=\sum_{P\in\mathscr{P}^n}\|\pi_{\mathscr{F}^n_P} z\|^2
\end{displaymath}
and hence for almost all $\omega\in\Omegait$
\begin{displaymath}
\|\pi_{\mathscr{P}^n} \xi^n(\omega)\|^2\le \sum_{P\in\mathscr{P}^n}C\cdot\ln|\mathscr{R}^n|=C\cdot|\mathscr{P}^n|\cdot\ln|\mathscr{R}^n|\:\:
\hbox{for eventually all}\:\: n\in\mathbb{M}.
\end{displaymath}    
This completes the proof.
\end{proof}


Let us finally illustrate the above concept in the classical case of gaussian white noise.

\begin{remark}\label{rem Noise Projection gwn}\rm
 Continuing from Remark \ref{white gaussian noise}, we  illustrate the behaviour of the lower bound for the constant $C$ in Proposition \ref{prop estimate projections} and Theorem  \ref{th Maximal equality} in the case of white gaussian noise and  polynomially growing number of fragments,
 i.e. $|\mathscr{R}^n|$ is asymptotically equivalent to $n^\kappa$. In this case the estimate for the norm of  noise projections takes the form  
\begin{eqnarray*}
    \|\pi_{\mathscr{P}^n} \xi^n(\omega)\|^2  &\leq &  \left( \frac{1}{\kappa} + 1\right)  \kappa 2 \sigma^2 D |\mathscr{P}^n| \ln n   = (1+\kappa) 2 \sigma^2 D      |\mathscr{P}^n|  \ln n, \\
   & & \mbox{ for almost each } \omega \mbox{ eventually.}     
\end{eqnarray*}
This underlines the dependency between the noise projections, the number of fragments, the noise variance, the dimension of the regression spaces and the size of the partitions.  
\end{remark}


%

\subsection{Discrete and Continuous Functionals}

We want to approximate functions $f$ on the continuous domain $S^\infty=[0,1)^d$ by estimates on  discrete finite grids $S^n$. The connections between the two settings are provided by the maps $\iota^n$ and $\delta^n$, introduced in \eqref{eq:delta} and \eqref{eq:iota}.
Note first that 
\begin{equation}\label{eq isometrie}
\langle\iota^n x,\iota^n y\rangle=\langle x,y\rangle/|S^n|\:\:\hbox{and}\:\: \|\iota^n x\|^2=\|x\|^2/|S^n|\quad \hbox{for}\:\: x,y\in\mathbb{R}^{S^n},
\end{equation}
where the inner product and norm on the respective left hand sides are those on $L^2(S^\infty)$ and on the right hand sides one has the Euclidean inner product and norm. Furthermore, one needs appropriate versions of the functionals (\ref{eq Potts functional}). Let now  $\mathfrak{S}^n$ be segmentation classes on the domains $S^n$ and $\mathfrak{S}\supset \iota^n\mathfrak{S}^n$ a segmentation class on $S^\infty$. 
Set
\begin{eqnarray*}
H_\gamma^n:\mathbb{R}^{S^n}\times\mathfrak{S}^n,\:\:H_\gamma^n(z,(\mathscr{P}^n,g_{\mathscr{P}^n}^n))&=&
\gamma|\mathscr{P}^n|+\|z-g_{\mathscr{P}^n}\|^2/|S^n|\\[1ex]
\tilde{H}^n_\gamma:L^2({S^\infty})\times\mathfrak{S},\:\:\tilde{H}^n_\gamma(f,(\mathscr{P},g_{\mathscr{P}}))&=&
\left\{
\begin{array}{cl}
\gamma|\mathscr{P}|+\|f-g_\mathscr{P}\|^2&\hbox{if}\:\:(\mathscr{P},g_\mathscr{P})\in\iota^n \mathfrak{S}^n,\\[1ex]
\infty                                                                     &\hbox{otherwise}.
\end{array}
\right.
\end{eqnarray*}
The two functionals are compatible. 
\begin{proposition} \label{lem:discrete_continuous_functionals} Let $n\in\mathbb{M}$ and $(\mathscr{P}^n,g_{\mathscr{P}^n}) \in\mathfrak{S}^n$ and $z^n\in\mathbb{R}^{S^n}$. Then
\begin{displaymath}\label{eq:functionals_equality} 
H_\gamma^n(z^n,(\mathscr{P}^n,g_{\mathscr{P}^n)}^n)=\tilde{H}_\gamma^n(\iota^n z^n,\iota^n(\mathscr{P}^n,g_{\mathscr{P}^n}^n)).
\end{displaymath}
If, moreover, $f\in L^2(S^\infty)$ then
\begin{displaymath}\label{eq argmin ist argmin}
(\mathscr{P}^n,g_{\mathscr{P}^n}^n)\in\argmin H^n_\gamma(\delta^n f,\cdot)
\:\hbox{if and only if}\:\iota^n(\mathscr{P}^n,g_{\mathscr{P}^n}^n)\in\argmin \tilde{H}^n_\gamma(f,\cdot)
\end{displaymath}
\end{proposition}

\begin{proof} The identity 
is an immediate 
consequence of (\ref{eq isometrie}). Hence let us turn to 
the equivalence of minimal points.  The key 
 is a suitable decomposition of the functional
$\tilde{H}^n_\gamma(f,\cdot)$. 
The map $\iota^n \delta^n$ is the
orthogonal projection of $L^2(S^\infty)$ onto the linear space $\mathscr{H}^n=\mathrm{span}\{\Einsbi_{I_{ij}:1\le i,j\le n}\}$,
and for any $(\mathscr{P},h)\in \iota^n \mathfrak{S}^n$ the function $h$ is in $\mathscr{H}^n$. Hence 
$$
    \|f-h\|^2 + \gamma|\mathscr{P}| = \|f - 
\iota^n  \delta^n f\|^2 + \|\iota^n  \delta^n f - 
h\|^2+ \gamma |\mathscr{P}|.
$$
The quantity  $\|f - \iota^n  \delta^n f\|^2$ does not depend on $(\mathscr{P},h)$. Therefore a pair 
$(\mathscr{P},h)$ minimises
$$
\|f - \iota^n  \delta^n f\|^2+ \|\iota^n  \delta^n f - 
h\|^2 + \gamma |\mathscr{P}|
$$
if and only if it minimises
$$
 \|\iota^n  \delta^n f - h\|^2 + \gamma |\mathscr{P}| = \tilde{H}_\gamma^n(\iota^n \delta^n 
f,\iota^n(\mathscr{P},h)).
$$
Setting $z ^n = \delta^n f$ in  (\ref{eq:functionals_equality}), this  completes the proof.
\end{proof}

\subsection{Upper Bound for Projective Segmentation Classes}\label{suse A Fundamental Inequality}

We compute an upper bound for the estimation error in a special setting: Choose in advance  a finite dimensional linear subspace $\mathscr{G}$ of $L^2(S^\infty)$. Discretization induces linear spaces $\delta^n\mathscr{G}=\{\delta^n f:f\in\mathscr{G}\}$ and   $\mathscr{G}^n_P=\{\Einsbi_P\cdot g:g\in \delta^n\mathscr{G}\}$, for any $P\subset S^n$, of functions on $S^n$. 
Let further for each $n \in \M$, a set $\mathscr{R}^n$ of admissible fragments and a family  $\mathfrak{P}^n$ of  partitions with fragments in $\mathscr{R}^n$  be given. Set $\mathfrak{G}^n := \left\{ \mathscr{G}_{\mathscr{P}} : \mathscr{P} \in \mathfrak{P}^n  \right\}$.
The induced segmentation class 
$$
\mathfrak{S}^n(\mathfrak{P}^n,\mathfrak{G}^n) = \left\{ (\mathfrak{P}^n,f) : \mathscr{P} \in \mathfrak{P}^n, f \in \mathscr{G}_{\mathscr{P}}\right\}
$$
will be called  \emph{projective ($\mathscr{G}$-) segmentation class} at stage $n$.


The following inequality is at the heart of later arguments since it controls the distance between the discrete $M$-estimates and the `true' signal.
\begin{lemma} \label{lem:fundamental_inequality}
Let for $n\in\mathbb{M}$ a $\mathscr{G}$-projective segmentation class $\mathfrak{S}^n$ over $S^n$ be given  and choose a signal $f\in L^2(S^\infty)$ and a vector $\xi^n\in \mathbb{R}^{S^n}$. Let further
\begin{displaymath}
(\hat{\mathscr{P}}^n,\hat{f}^n)\in\argmin_{(\mathscr{Q},h)\in\mathfrak{S}^n} H_\gamma^n(\delta^n f+\xi^n,(\mathscr{Q},h))
\end{displaymath}
and $(\mathscr{Q},h)\in\mathfrak{S}^n$. Then
\begin{equation}\label{eq fundamental inequality}
\|\iota^n \hat{f}^n-f\|^2\le 2\gamma(|\mathscr{Q}|-|\hat{\mathscr{P}}^n|)+3\|\iota^n h-f\|^2+ \frac{16}{n^d}
 \left( \|\pi_{\hat{\mathscr{P}}^n  }\xi^n\|^2 + \|\pi_{\mathscr{Q}  }\xi^n\|^2 \right).
\end{equation}
\end{lemma}
\begin{proof}
Since $(\hat{\mathscr{P}}^n,\hat{f}^n)$ is a minimal point of  $H_\gamma^n(\delta^n f+\xi^n,\cdot))$ the embedded segmentation $\iota^n(\hat{\mathscr{P}}^n,\hat{f}^n)$ is a minimal point of $\tilde{H}_\gamma^n( f+\iota^n\xi^n,\cdot))$ by Proposition \ref{lem:discrete_continuous_functionals} and hence
\begin{displaymath}
\gamma|\hat{\mathscr{P}}^n|+\|(\iota^n \hat{f}^n-f)-\iota^n\xi^n\|^2\le \gamma|\mathscr{Q}|+\|(\iota^n h-f)-\iota^n\xi^n\|^2.
\end{displaymath}
 Expansion of squares yields that
 \begin{eqnarray*}
&&\gamma|\hat{\mathscr{P}}^n|+\|\iota^n \hat{f}^n-f\|^2+2\langle\iota^n \hat{f}^n-f,\iota^n\xi^n\rangle+\|\iota^n\xi^n\|^2\\
&\le& \gamma|\mathscr{Q}|+\|\iota^n h-f\|^2+2\langle\iota^n h-f,\iota^n\xi^n\rangle+\|\iota^n\xi^n\|^2
\end{eqnarray*}
and hence
\begin{equation} \label{eq:expansion}
\|\iota^n \hat{f}^n-f\|^2\le \gamma(|\mathscr{Q}|-|\hat{\mathscr{P}}^n|)+\|\iota^n h-f\|^2+2\langle\iota^n h-\iota^n \hat{f}^n,\iota^n \xi^n\rangle.
\end{equation}
By definition $h \in \mathscr{F}_{\mathscr{Q}}$ and $\hat{f}^n \in \mathscr{F}_{\hat{\mathscr{P}}^n}$ which implies that 
$h-\hat{f}^n\in\mathscr{F}^\prime = \mathrm{span} ( \hat{\mathscr{P}}^n, \mathscr{F}_{\mathscr{Q}})
$
and hence 
$\pi_{\mathscr{F}^\prime}(\hat{f}^n -h)=\hat{f}^n -h$.  
We proceed with
\begin{eqnarray*}
&&|\langle\iota^n h-\iota^n \hat{f}^n,\iota^n\xi^n\rangle|=|S^n|^{-1}|\langle\pi_{\mathscr{F}^\prime }(\hat{f}^n -h),\xi^n\rangle|
=|S^n|^{-1}|\langle  h-\hat{f}^n,\pi_{\mathscr{F}^\prime }\xi^n\rangle|\\
&\le&\|\iota^n \hat{f}^n-\iota^n h\|\cdot|S^n|^{-1/2}
\cdot \|\pi_{\mathscr{F}^\prime}\xi^n\|\\
&\le&|S^n|^{-1/2}\|\pi_{\mathscr{F}^\prime}\xi^n\|\cdot\|\iota^n \hat{f}^n-f\|   +   |S^n|^{-1/2}\|\pi_{\mathscr{F}^\prime}\xi^n\|\cdot
\|f-\iota^n h\|.
\end{eqnarray*}
Since $ab\le a^2+b^2/4$,  we conclude
\begin{eqnarray*}
|\langle\iota^n h-\iota^n \hat{f}^n,\iota^n\xi^n\rangle| &\le& \|\iota^n \hat{f}^n-\iota^n h\|^2\//4+\|f-\iota^n h\|^2/4+
2\|\pi_{\mathscr{F}^\prime}\xi^n\|^2/|S^n| \\
&\leq& \|\iota^n \hat{f}^n-\iota^n h\|^2\//4+\|f-\iota^n h\|^2/4+
4 \left( \|\pi_{\hat{\mathscr{P}}^n} \xi^n\|^2 +  \|\pi_{\mathscr{Q}}\xi^n\|^2\right)/|S^n|
\end{eqnarray*}
Putting this into inequality~(\ref{eq:expansion}) results in
\begin{eqnarray*}
\|\iota^n \hat{f}^n-f\|^2&\le& \gamma(|\mathscr{Q}|-|\hat{\mathscr{P}}^n|)+\|\iota^n h-f\|^2+\|\iota^n \hat{f}^n-f\|^2/2+\|f-\iota^n h\|^2/2\\
&+&8
\left( \|\pi_{\hat{\mathscr{P}}^n }\xi^n\|^2 + \|\pi_{\mathscr{Q} }\xi^n\|^2 \right)            
/ |S^n|,
\end{eqnarray*}
which implies the asserted inequality.
\end{proof}

\section{Consistency} \label{se Consistency}
In this section we complete the abstract considerations and summarize the preliminary work in two theorems on consistency. The first one  concerns the desired $L^2$-convergence of estimates to the `truth', and the second one provides convergence rates.
\subsection{$L^2$-Convergence}
We will prove now that the estimates of images converge almost surely to the underlying true signal in $L^2(S^\infty)$ for almost all observations. We adopt the projective setting introduced in Section \ref{suse A Fundamental Inequality}.
Let us make some agreements in advance.
\begin{hypothesis} \label{hyp rate assumptions}
Assume that
\begin{itemize}
\item[(H\ref{hyp rate assumptions}.1)] there are $\kappa>0$ and $C>0$ such that $|\mathscr{R}^n|\ge C\cdot n^\kappa$ eventually,
\item[(H\ref{hyp rate assumptions}.2)] there is a real number $\beta >0$ such that, 
for each $n\in\mathbb{M}$ and real numbers $\mu_s$, $s\in S^n$, and each $c\in\mathbb{R}_+$, the inequality
\begin{equation*}
\mathbb{P}\left(\left|\sum_{s\in S_n}\mu_s\xi_s^n\right|\ge c\right)\le 2\cdot \exp\left(-\displaystyle \frac{c^2}{\beta\sum_{s \in S_n}\mu_s^2}\right)
\end{equation*}
holds,
\item[(H\ref{hyp rate assumptions}.3)] the positive sequence $(\gamma_n)_{n\in\mathbb{N}}$ satisfies 
\begin{equation*}
\gamma_n\rightarrow0 \mbox{ and } \gamma_n > C\cdot \frac{\ln |\mathscr{R}^n|}{|S^n|}, \mbox{ for eventually all } n 
\end{equation*}
with $C=\beta D (\kappa+1)/ \kappa  $,  and $D$ is, like in Proposition \ref{prop estimate projections}, an upper bound for the dimension of the linear spaces $\mathscr{F}_P$.
\end{itemize}
\end{hypothesis}

\begin{bemerkung}
Note that the condition $\gamma_n\cdot|S^n|/\ln n\rightarrow\infty$ implies the second part of (H\ref{hyp rate assumptions}.3) by (H\ref{hyp rate assumptions}.1). It was used for example in \fcite{Friedrich:2005} or \fcite{LiebscherBoysenMunk:2009,BoysenLiebscherMunkWittich:2007}.
\end{bemerkung}

Given a signal $f\in L^2(S^\infty)$ we must assure that our setting actually allows for good approximations of $f$ at all. If so, least squares estimates
are consistent.
\begin{theorem}\label{thm consistency: L2-convergence }
Assume that Hypothesis \ref{hyp rate assumptions} holds.
Let $f\in L^2(S^\infty)$ and suppose 
\begin{equation} \label{eq:approx_condition}
\lim_{k\rightarrow \infty}\limsup_{n\rightarrow\infty}\inf_{(\mathscr{Q},h)\in\mathfrak{S}^n, |\mathscr{Q}|\le k}\|\iota^n h-f\|^2=0.
\end{equation}
Then
\begin{displaymath}
\|\iota^n \hat{f}^n(\omega)-f\|^2\longrightarrow 0\quad\hbox{as}\:n\rightarrow\infty\hbox{  for almost all  }\omega \in\Omegait.
\end{displaymath}
\end{theorem}
We formulate part of the proof separately, since it will be needed later once more.
\begin{lemma}\label{la first inequality}
We maintain the assumptions of  Theorem \ref{thm consistency: L2-convergence }.  
Then, given $k > 0$,
\begin{equation}\label{eq first inequality}
\|\iota^n \hat{f}^n(\omega)-f\|^2\le 3\gamma_n\cdot k+3\|\iota^n h-f\|^2\hbox{for all $(\mathcal{Q},h) \in  \mathfrak{S}^n$ such that  $|\mathcal{Q}| \leq k$ }
\end{equation}
eventually for all  $n \in \mathbb{N}$ and for almost all $\omega\in\Omegait$.
\end{lemma}
\begin{proof} Lemma~\ref{lem:fundamental_inequality} yields
$$
    \|\iota^n \hat{f}^n (\omega) -f\|^2 \leq 2 \gamma_n \left( |\mathscr{Q}| - |\mathscr{P}^n|\right) + 3 \|\iota^n h -f \|^2 + \frac{16}{n^d} \left( \|\pi_{\hat{\mathscr{P}}^n} \xi \|^2  +\|\pi_{\mathscr{Q}} \xi \|^2 \right) 
$$
and application of Proposition \ref{prop estimate projections} implies that for any
real number $C^\prime>\frac{\kappa+1}{\kappa} \beta D$, the following inequality holds for  almost all $\omega\in \Omegait$ 
%
%
%
%
%
%
\begin{eqnarray*}
\|\iota^n \hat{f}^n (\omega) -f\|^2 &\leq &2 \gamma_n k + 3 \|\iota^n h -f\|^2 + 16 C^\prime \left( \frac{\ln (|\mathscr{R}^n|)}{n^d}\right) \cdot 
\left( 
|\mathscr{Q} |+ |\hat{\mathscr{P}}^n| \right) - 2 \gamma_n \cdot |\hat{\mathscr{P}}^n| \\
&\leq& 2 \gamma_n k + 3 \|\iota^n h -f\|^2 + 16C^\prime \frac{\ln |\mathscr{R}^n|}{n^d} k + |\mathscr{P}^n| \left(  8 C^\prime \frac{\ln |\mathscr{R}^n|}{n^d}  - 2 \gamma_n \right)
\end{eqnarray*}
For $\gamma_n $ satisfying  Hypothesis (H\ref{hyp rate assumptions}.3), the term in parenthesis is negative.
Therefore (\ref{eq first inequality}) holds and the assertion is proved.
\end{proof}



Theorem \ref{thm consistency: L2-convergence } follows now easily.
\begin{proof}[ Proof of Theorem \ref{thm consistency: L2-convergence }]
The following formulae hold almost surely. Lemma \ref{la first inequality} implies that, for 
$$
    \| \iota^n \hat{f}^n -f \|^2 \leq  3 \gamma_n \cdot k + 3 \cdot \inf_{(\mathscr{Q},h)\in\mathfrak{S}^n,|\mathscr{Q}|\le k} 
    \left( \|\iota^n h -f\|^2 \right) \hbox{ eventually}
$$
Therefore  
\begin{eqnarray*}
   \limsup_{n\rightarrow\infty} \| \iota^n \hat{f}^n -f \|^2  &\leq  & \limsup_{n\rightarrow\infty} \left( 3 \gamma_n \cdot k + 3 \cdot \inf_{(\mathscr{Q},h)\in\mathfrak{S}^n,|\mathscr{Q}|\le k} 
    \left( \|\iota^n h -f\|^2 \right) \right)  \\
    & = & 0  \; + \; 3 \cdot \limsup_{n\rightarrow\infty} \inf_{(\mathscr{Q},h)\in\mathfrak{S}^n,|\mathscr{Q}|\le k} 
    \left( \|\iota^n h -f\|^2 \right)
\end{eqnarray*}
By assumption  (\ref{eq:approx_condition}), the right hand side converges to 0 as $k$ tends to $\infty$.  Hence
$$
    \limsup_{n \rightarrow \infty} \|\iota^n \hat{f}^n-f\|^2 = 0,
$$
which completes the proof. 
\end{proof}
\subsection{Convergence Rates}
The final abstract result  provides  almost sure 
convergence rates in the general setting. 
\begin{theorem} \label{th convergence rates}
Suppose that Hypothesis \ref{hyp rate assumptions} holds 
and assume further that there are real numbers $\alpha, C >0,\varrho\ge0 $, and a sequence $(F_n)_{n\in \mathbb{N}}$ with $\lim_{n\rightarrow\infty} F_n=\infty$ such that
\begin{eqnarray}\label{eq Bedingung mit F}
&&\|\iota^n h-f\| \le C\cdot\left(\frac{k^\varrho}{F_n}+\frac{1}{k^\alpha}\right)\\
&&\hbox{for all } n\in\mathbb{M}\hbox{ and } k,\hbox{ and some }(\mathscr{Q},h)\in\mathfrak{S}^n\:\:\hbox{with } |\mathscr{Q}|\le k.\nonumber
\end{eqnarray}
Then
\begin{equation} \label{eq:minimax_discretisation_estimate}
\|\iota^n \hat{f}^n(\omega)-f\|^2=O\left( \gamma_n^{\frac{2\alpha}{2\alpha +1}} \right) + O\left(F_n^{-\frac{2\alpha}{\alpha+\varrho}}\right)
\hbox{  for almost all }\omega \in\Omegait.
\end{equation}
\end{theorem}
%
\begin{proof} 
Let $(k_n)_{n\in\M}$ be a sequence in $\mathbb{R}_+$. Recall from  Lemma \ref{la first inequality} that 
$$
    \|\iota^n \hat{f}^{n} - f\|^2 \le 2 \gamma_n \cdot k_n + 3 \cdot \|\iota^n  h  - f \|_2^2
$$
for  sufficiently large $n \in \M$ and any $(\mathscr{Q},h) \in \mathfrak{S}^n$ with $|\mathscr{Q}|\le k_n$ on a set of $\omega$ of full measure. 
The following arguments hold for all such $\omega$. We will write $C$ for constants; hence the $C$
 below may differ. 
 
\noindent Since $(a+b)^2\le 2(a^2+b^2)$, assumption  (\ref{eq Bedingung mit F}) implies that
\begin{equation} \label{eq:error_decomposition}
    \|\iota^n \hat{f}^{n} - f\|^2 \leq C\left( \gamma_n \cdot k_n+  \frac{k_n^{2 \varrho}}{F_n^2} + 
    \frac{1}{k_n^{2 \alpha}} \right).
\end{equation}
This decomposition of the error can be interpreted as follows: the first term corresponds to an estimate of the error due to the noise, the second term corresponds to the discretization while the third term can be directly related to the approximation error of the underlying scheme, in the continuous domain. 

One has free choice of the parameters  $k_n$.
We enforce the same decay rate for the first and third term  setting
$
\gamma_n k_n = k_n^{-2 \alpha}
$
Then, in view of (\ref{eq:error_decomposition}),
\begin{equation} \label{eq:minimax}
    \|\iota^n \hat{f}^{n} - f\|^2 \leq C\left(\gamma_n^{\frac{2 \alpha}{2 \alpha +1}}+  \frac{\gamma_n^{- \frac{2 \varrho}{2 \alpha + 1}}}{F_n^2}\right).
\end{equation}
To get the same rate for the discretisation and the approximation error set
\begin{displaymath} \label{eq:relation_kn_Fn}
    \frac{k_n^{2 \varrho}}{F_n^2} = \frac{1}{k_n^{2 \alpha}} \mbox{ or equivalently } k_n = F_n^{\frac{1}{\varrho + \alpha}},
\end{displaymath}
which, together with  estimate (\ref{eq:error_decomposition}),  yields
\begin{equation} \label{eq:discretisation}
    \|\iota^n \hat{f}^{n} - f\|^2 \leq C\left( \gamma_n F_n^{\frac{1}{\varrho + \alpha}}  +F_n^{-\frac{2 \alpha}{ \alpha + \varrho}} \right). 
\end{equation}
Straightforward calculation gives
\begin{displaymath}
\gamma_n ^{\frac{2\alpha}{2\alpha+1}} \ge  \frac{\gamma_n^{-\frac{2\varrho}{2\alpha+1}}}{F_n^2}
\hbox{ if and only if }\gamma_n F_n^{\frac{1}{\alpha + \varrho}} \ge \frac{1}{F_n^{\frac{2 \alpha}{\alpha + \varrho}}}
\end{displaymath}
Hence, the first term on the right hand side of inequality  (\ref{eq:minimax}) dominates the second one if and only this holds in  inequality  (\ref{eq:discretisation}). We discriminate between the two  cases $\ge$ and $<$. The first one is 
\begin{equation}\label{eq 1Fall}
\gamma_n ^{\frac{2\alpha}{2\alpha+1}}\ge\frac{\gamma_n^{-\frac{2\varrho}{2\alpha+1}}}{F_n^2}.
\end{equation}
Combination with (\ref{eq:minimax}) results in 
\begin{equation} \label{estimate 1Fall}
    \|\iota^n \hat{f}^{n} - f\|_2^2\le C\cdot\gamma_n^{\frac{2\alpha}{2\alpha+1}}
\end{equation}
for some $C>0$. In view of the equivalence, replacement of $\ge$ by $<$ in (\ref{eq 1Fall}), results in
\begin{displaymath} \label{eq 2Fall}
\gamma_n F_n^{\frac{1}{\alpha+\varrho}} < F_n^{-\frac{2\alpha}{\alpha+\varrho}}.
\end{displaymath}
which, together with  
estimate (\ref{eq:discretisation}),  gives for some $C>0$ that
  \begin{equation} \label{estimate 2Fall}
    \|\iota^n \hat{f}^{n} - f\|^2 \leq C\cdot  F_n^{-\frac{2 \alpha}{\alpha + \varrho}}.
  \end{equation}
 Combination of  (\ref{estimate 2Fall}) and (\ref{estimate 1Fall})  completes the proof of (\ref{eq:minimax_discretisation_estimate}).
\end{proof}

\begin{remark}\rm Let us continue from Remark \ref{rem Noise Projection gwn}. If  $|\mathscr{R}^n| \sim n^\kappa$ and noise is white Gaussian  with  $\beta = 2 \sigma^2$ then
 Hypothesis (H\ref{hyp rate assumptions}.3)
boils down to
$$
\gamma_n\longrightarrow 0\hbox{  and  }    \gamma_n >  2 (\kappa+1) \sigma^2 D\cdot \frac{\ln n}{n^d}.
$$
Setting $\varepsilon_n = \sigma/n^{d/2}$,  the estimate (\ref{eq:minimax_discretisation_estimate}) then reads
$$
    \|\iota^n \hat{f}^n (\omega) -f \|^2 = O\left(   \left( \varepsilon_n^2 \left| \ln  \varepsilon_n \right| \right)^{\frac{2 \alpha}{2 \alpha + 1}} \right),
$$ 
as long as the growth of $F_n$ is sufficient. This is strongly connected with the optimal minimax rates from model selection, 
which bound the expectations of the left hand side, see for instance \fcite{BirgeMassart:1997}. 
\end{remark}
%

\section{Special Segmentations}\label{se Special Segmentations}

We are going now to exemplify  the abstract Theorem  \ref{th convergence rates} by way of  typical partitions and spaces of functions. On the one hand, this extends a couple of already existing results and, on the other hand, it  illustrates the wide range of possible applications.

\subsection{One Dimensional Signals - Interval Partitions}
\label{sec:1D_interval_partitions}

Albeit focus of this paper is on two  or more dimensions, we start with one dimension.
There are at least two reasons for that: illustration of  the abstract results by choices of the (seemingly) most elementary example, and  to generalize results like some of those  in \fcite{LiebscherBoysenMunk:2009,BoysenLiebscherMunkWittich:2007}  to classes of  piecewise Sobolev functions.

To be definite,  let   $S^n=\{1,\,\ldots\,,n\}$  and let $\mathscr{R}^n = \{[i,j] :1\le  i\le j \le n\}$ be the  discrete intervals of admissible fragments. Then $\mathfrak{P}^n$   is  the collection of partitions of  $S^n$ into intervals.
Plainly, $|\mathscr{R}^n|=(n+1)n/2$ and 
$ |\mathfrak{P}^n| =2^{n-1}$. We deal with approximation by polynomials. To this end and in accordance with Section \ref{suse A Fundamental Inequality}, we  choose the finite dimensional linear subspace $\mathscr{F}_p \subset L^2([0,1))$ of polynomials of maximal degree $p$. The induced segmentation  classes
$\mathfrak{S}^n(\mathfrak{P}^n,\mathfrak{F}^n)$ consist of piecewise polynomial functions relative to partitions in $\mathfrak{P}^n$.  
 
The signals to be estimated will be members of 
the  fractional Sobolev space $W^{\alpha,2}((0,1))$ of order $\alpha>0$.   The main task is to verify Condition (\ref{eq Bedingung mit F}).  Note that this class of functions is slightly larger than the classical H\"{o}lder spaces of order $\alpha$ usually treated.
For results in the case of equidistant partitioning, we refer, for instance, to  \fcite{gyoerfi2002distribution} Section 11.2.

For the  following lemma, we adopt classical arguments from approximation theory.
\begin{lemma}  
\label{la:piecewise_polynomial_approximation_discretisation} 
For any $f \in W^{\alpha,2}((0,1))$,  with 
$p<\alpha<p+1$, there is $C>0$ such that for all $k \leq n \in \N$, 
there is
 $(\mathscr{P}^n_k,h^n_k,) \in \mathfrak{S}^n$, such that 
$|\mathscr{P}^n_k| \leq k$ and which satisfies
\begin{equation}
\label{eq:approximation_discretisation}
\|f-\iota^n h^n_k\| \leq C \cdot \left( \frac{1}{k^\alpha} +    
\frac{k}{n} \right)
\end{equation}
\end{lemma}

For the proof, let us introduce partitions 
$\mathscr{I}_k=\{[(i-1)/k,i/k) : i=1,\cdots,k\}$ of $[0,1)$ into $k$ intervals, each of length $1/k$.
\begin{proof} Let $f \in W^{\alpha,2}((0,1))$.
From classical approximation theory (see e.g. \cite{DeVoreLorentz:1993}, 
Chapter 12, Thm. 2.4), we learn  that there is  $C>0$ such that
there is a piecewise polynomial function $h_k$ 
of degree at most $p$ such that
$$
\|f - h_k\| \leq \frac{C}{k^\alpha}.
$$
For each $i=1,\ldots,k$,  let $h_{k,i}$ denote the restriction of $h_k$ to $I_i=((i-1)/k,i/k)$. We 
consult the Bramble-Hilbert lemma (for a version corresponding to our 
needs, we refer to Thm. 6.1 in \cite{DupontScott:1980})
and find   $C>0$, such that
$$
    |f-h_{k,i}|_{W^{1,2}(I_i)} \leq C \cdot |f|_{W^{1,2}(I_i)}  \:\hbox{ for each  }i=1,\ldots,k.
$$
This yields for some $C>0$, independent of  $k$ and $n$, that
$$
   |h_{k,i}|_{W^{1,2}(I_i)} \leq  |f - h_{k,i}|_{W^{1,2}(I_i)} + 
|f|_{W^{1,2}(I_i)}  \leq C\cdot |f|_{W^{1,2}(I_i)} \mbox{ for all  } i=1,\cdots,k.
$$
We turn  now to the  piecewise 
constant  approximation on the  partition  $\mathscr{I}_n$.
We split $[0,1)$ into the  union $J_k^n$ of those intervals in $\mathscr{I}_n$ which do not 
contain knots  $i/k$ and the union $K_k^n$ of those intervals in $\mathscr{I}_n$ which  do
contain knots  $i/k$.  For $I \in \mathscr{I}_k$ and $I \subset J_k^n$, we have 
\begin{displaymath}
|h_{k,i}|_{W^{1,2}(I)} \leq  C |f|_{W^{1,2}(I)} |\;\;\hbox{  if and only if      }\;\;
 |h^\prime_{k,i}|^2_{L^2(I)} \leq  C^2 \cdot |f^\prime|^2_{L^2(I)}. 
\end{displaymath}
This implies
$$
    \sum_{I \subset J_n^k} |h_{k,i}^\prime|^2_{L^2(I)} \leq C^2 \sum_{I \subset J_n^k} |f^\prime|^2_{L^2(I)} \leq C^2 |f^\prime|_{L^2([0,1])},
$$
which in turn leads to
$$
    |h_k|_{W^{1,2}(J_n^k)} \leq C^2 |f|_{W^{1,2}((0,1))}.
$$
Hence we are ready  to conclude that for some constant $C>0$,
\begin{equation} \label{eq:1D_Jnk}
   \|h_k- \iota^n \delta^n h_k\|_{L^2(J_k^n)} \leq C/n.
\end{equation}
For $I \in \mathscr{I}_k$ and $I \subset K_k^n$, we use the fact that $h_k^n \leq 2C \cdot \|f\|_{L^\infty}([0,1])$ and deduce 
$$
   \|h_k-\iota^n \delta^n h_k\|_{L^2(I)} \leq 2 C \|f\|_{L^\infty(I)}/n.
$$
Summation over all intervals included in $K_k^n$ results in
$$
    \|h_k- \iota^n \delta^n h_k\|_{L^2(K_k^n)}  \leq  C \cdot 
k/n.
$$
\noindent This yields for the entire interval $[0,1)$ that
$$
   \|f- \iota^n \delta^n h_k \| \leq  \|f-h_k\| + \|h_k - \iota^n \delta^n h_k\| \leq C \left( \frac{k}{n} +  
\frac{1}{k^\alpha}\right).
$$
With  $h_k^n = \delta^n h_k$, this completes the proof.
\end{proof}


Piecewise smooth functions have only a very low Sobolev regularity. Indeed, recall that piecewise smooth functions belong to $W^{\alpha,2}((0,1))$ only for $\alpha > 1/2$. In order to overcome this limitation, we consider a larger class of functions, the class
of piecewise Sobolev functions.
\begin{definition} \label{def:piecewise_Sobolev} Let $\alpha>1/2$ be a  real number, $J \in \N$, and $x_0 =0 <x_1 < \cdots < x_{J+1} = 1$. A function $f$ is said to be \emph{piecewise $W^{\alpha,2}([0,1])$ with $J$ jumps}, relative to the partition  $\{[x_i,x_{i+1}): i=1,\cdots,J\}$ if 
$$
  f |_{(x_i,x_{i+1})}  \in W^{\alpha,2}\left((x_i,x_{i+1}) \right)
$$ 
\end{definition}
\begin{bemerkung} Definition \ref{def:piecewise_Sobolev} is consistent, due to the Sobolev embedding theorem. For  an open interval $I$   of $\R$, $ W^{\alpha,2}(I) $ is continuously embedded into 
$ \mathscr{C}\left(I^a\right)$, the space of  uniformly continuous functions on the closure
  $I^a$ of $I$.
\end{bemerkung}

We conclude from  Lemma \ref{la:piecewise_polynomial_approximation_discretisation}:
\begin{lemma} \label{la:piecewise_polynomial_approximation_discretisation2}
Let $f$ be piecewise-$ W^{\alpha,2}([0,1))$ with 
$J$ jumps and   with $p<\alpha<p+1$. Then  there are  $C>0$ and 
$(\mathscr{P}^n_k,h^n_k) \in \mathfrak{S}^n$, such that $|\mathscr{P}^n| 
\leq k$ and
\begin{equation}
\label{eq:approximation_discretisation}
\|f-h^n_k\| \leq C \cdot \left( \frac{1}{k^\alpha} +    \frac{k}{n} + 
\frac{J}{n} \right).
\end{equation}
\end{lemma}

\begin{proof} With the same arguments as in the previous proof we just 
have to include the error made at each jump of the original piecewise regular function. 
More precisely, we use a similar splitting into $J_k^n$ and $K_k^n$ where $K_k^n$ also
contains the intervals containing $x_i$ for $i=1,\cdots,J$. Since there are at most $k+J$ intervals in $K_k^n$, this gives estimate
\eqref{eq:approximation_discretisation}.
\end{proof}

By Lemma \ref{la:piecewise_polynomial_approximation_discretisation2}, 
a piecewise Sobolev function satisfies  Condition \eqref{eq Bedingung mit F} with $\rho = 1$ and $F_n = n$ and therefore 
Theorem \ref{th convergence rates} applies. In summary

\begin{theorem} \label{th:1d_consistency}  Let $f$ be a  piecewise $ W^{\alpha,2}([0,1])$ function, such that $0<\alpha<p+1$, where $p$ is the maximal degree of the approximating polynomials. We assume further that (H\ref{hyp rate assumptions}.3) holds
and that the noise variables $\xi_s^n$ from Section \ref{suse The Regression Model} satisfy (\ref{eq MainTailEstimate}). Then
\begin{equation} \label{eq:1D_consistency}
\|\iota^n \hat{f}^n(\omega)-f\|^2=O\left( \gamma_n^{\frac{2\alpha}{2\alpha +1}}\right), 
\hbox{  for almost all }\omega \in\Omegait.
\end{equation}
\end{theorem}

\begin{proof} Let us  check  the assumption in Theorem \ref{th convergence rates}. 
Since $|\mathscr{R^n}|=(n-1)n/2$, Hypothesis  {\em (H\ref{hyp rate assumptions}.1)} holds with $\kappa = 2$. Hypothesis  {\em (H\ref{hyp rate assumptions}.2)}
and  {\em (H\ref{hyp rate assumptions}.3)} were required separately. Finally, Condition (\ref{eq Bedingung mit F}) holds with $\varrho=1$ and $F_n=n$ by Lemma  \ref{la:piecewise_polynomial_approximation_discretisation2}. Finally, Hypothesis  {\em (H\ref{hyp rate assumptions}.3)} completes the proof.
\end{proof}

Let $\mathscr{C}^1([0,1])$ denote the set of continuously differentiable functions.
For $p\in \N$, $\alpha \in (p,p+1]$, a function $f \in \mathscr{C}^p([0,1])$  is said to be $\alpha$\emph{-Hölder} if there is  $C >0$ such that
$$
     |f^{(p)}(x)-f^{(p)}(y)|  \leq C |x-y|^{\alpha-p}\; \mbox{ for any } \; x,y \in [0,1], x \neq y. 
$$
The linear space of $\alpha$-H\"older functions will be  denoted by $\mathscr{C}^{\alpha}([0,1])$ if $\alpha \in \N$ and $\mathscr{C}^{\alpha-1,1}([0,1])$ if $\alpha \in \N$.

\begin{bemerkung} \label{rem:simplified_estimate} Choose $\gamma_n = C \ln n/n$ with large enough $C$,  independently of $f$. Then the almost sure estimates~\eqref{eq:1D_consistency} of the estimation error simplifies to 
\begin{equation} \label{eq:simplified_estimate}
\|\iota^n \hat{f}^n(\omega)-f\|^2 = O\left( \frac{\ln n}{n} \right)^{\frac{\alpha}{2\alpha+1}} \hbox{  for almost all }\omega \in\Omegait.
\end{equation}
   These convergence rates are, up to the logarithmic factor, the optimal
rates for mean square error  in the Hölder classes $\mathscr{C}^\alpha([0,1])$. 
Thus, our  estimate adapts automatically  to the smoothness  of the signal. 

\end{bemerkung}

\subsection{Wedgelet Partitions}


Wedgelet decompositions are content-adapted partitioning methods based on elementary geometric atoms, called \emph{wedgelets}. A wedge results from the splitting of  a square into two pieces by a straight line and in our setting a wedgelet will be a piecewise polynomial function over a wedge partition. 
The discrete setting requires a careful treatment.  We adopt the discretization scheme from \fcite{Friedrich:2007}, which relies on the digitalization  of lines from \fcite{Bresenham:1965}.  
This discretization   differs from that in  \fcite{Donoho:1999}, where all pairs of  pixels on the boundary of a discrete square are used as endpoints of line segments.  One of the main reasons for our special choice  is  an efficient  algorithm which returns exact solutions of the functional (\ref{eq Potts functional}). It relies on  rapid moment computation, based on lookup tables, cf. \fcite{Friedrich:2007}.

\subsubsection*{Wedgelet partitions}

Let us first recall the relevant concepts and definitions. Only the case of dyadic wedgelet partitions will be  discussed. Generalisations  are straightforward but technical.

We start from discrete dyadic squares $S^{m}=\{1,\,\ldots\,,m\}^2$ with $m\in\M=\{2^p:p\in\N_0\}$. \emph{Admissible fragments} are dyadic squares of the form
\begin{displaymath}
[(i-1)\cdot 2^{q},i\cdot 2^{q})\times[(j-1)\cdot 2^{q},j\cdot 2^{q}),\quad 1\le i,j\le 2^{p-q}, \:0\le q\le p. 
\end{displaymath}
The collection  of dyadic squares  can be interpreted as the set of leaves of a quadtree where each internal node has exactly four children obtained by subdividing one square into four.


Digital lines in $\Z^2$ are defined for  angles $\vartheta\in (-\pi/4,3\pi/4]$.  Let
\begin{displaymath}
d(\vartheta)=\max\{|\cos \vartheta|,|\sin\vartheta|\},\quad v(\vartheta)=
\left\{
\begin{array}{ll}
(-\sin\vartheta,\cos\vartheta)&\hbox{if  }|\cos\vartheta|\ge |\sin\vartheta|\\
(\sin\vartheta,-\cos\vartheta)&\hbox{otherwise}
\end{array}
\right..
\end{displaymath}
The \emph{digital line through the origin in direction} $\vartheta$ is defined as
\begin{displaymath}
L^0_\vartheta=\{ s\in \Z^2:-d(\vartheta)/2<\langle s,v(\vartheta)\rangle\le d(\vartheta)/2  \}.
\end{displaymath}
Lines parallel to $L^0_\vartheta$ are  shifted versions
\begin{displaymath}
L^r_\vartheta=\{ s\in \Z^2:(r-1/2)d(\vartheta)<\langle s,v(\vartheta)\rangle\le (r+1/2)d(\vartheta)  \}
\end{displaymath}
with the \emph{line numbers} $r\in \Z$.  One  distinguishes between  \emph{flat} lines where  $\cos\vartheta\ge\sin\vartheta$ and  \emph{steep} lines where  $\cos \vartheta <\sin\vartheta$. 
For $x\in\R$, set $\mathtt{round}(x)=\max\{i\in\Z:i\le x+1/2\}$, let $y_\vartheta(x)=\mathtt{round}(x\cdot\tan\vartheta)$ and 
$x_\vartheta(x)=\mathtt{round}(y\cdot\cot\vartheta)$. According to Lemma 2.7 in \fcite{Friedrich:2007},
\begin{eqnarray*}
L^r_\vartheta&=&(0,r)+\{(x,y_\vartheta(x):x\in\Z)\}\hbox{  for flat lines},\\
L^r_\vartheta&=&(r,0)+\{(x_\vartheta(y),y:y\in\Z)\}\hbox{  for steep lines}.
\end{eqnarray*}
By Lemma 2.8 in the same reference, all parallel lines partition $\Z^2$. We are now ready to define wedgelets. Let $Q$ be a square in $\Z^2$ and $L_\vartheta^r$ a line with $L_\vartheta^r\cap Q\ne\emptyset$ and $L_\vartheta^{r+1}\cap Q\ne\emptyset$. A \emph{wedge split} is a partition of $Q$ into the \emph{lower} and \emph{upper wedge}, respectively, given by
\begin{equation} \label{eq:wedge_splitting}
W^{r,l}_{\vartheta}=\bigcup_{k\le r}L_\vartheta^k\cap Q,\quad W^{r,u}_{\vartheta}=\bigcup_{k> r}L_\vartheta^k\cap Q. 
\end{equation}

Let $\mathscr{Q}$ be a partition of some domain $S^m$ into squares. Then a  \emph{wedge partition}  of $S^m$ is obtained  replacing some of these squares by the two wedges of a wedge split. It is called \emph{dyadic} if  $m\in\M$, and the squares $Q\in\mathscr{Q}$ are dyadic.

We assume  that a finite set $\Thetait$ of angles is given. 
 The set $\mathscr{R}^m$ of admissible segments consists of wedges obtained by  wedge splits of  dyadic squares, given by (\ref{eq:wedge_splitting}) and  for $\theta \in \Thetait$, or by dyadic squares.

Focus  is on piecewise polynomial approximation of low order.
The induced segmentation  classes
$\mathfrak{S}^m$ consist of piecewise polynomial functions relative to a wedgelet partition. 
The cases of piecewise constant (original wedgelets) and piecewise linear polynomials (platelets) will be treated explicitly.


\begin{figure}[t!!]
\begin{tabular}{cc}
\includegraphics[width=5.8cm]{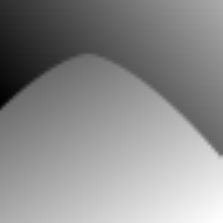}
\quad &
\includegraphics[width=5.8cm]{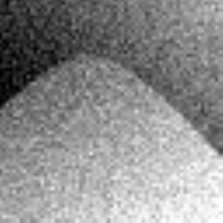}
\\
\includegraphics[width=5.8cm]{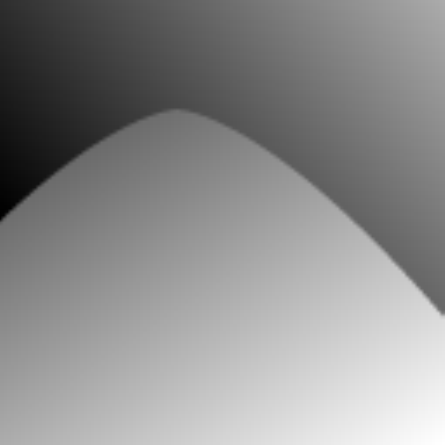}
\quad &
\includegraphics[width=5.8cm]{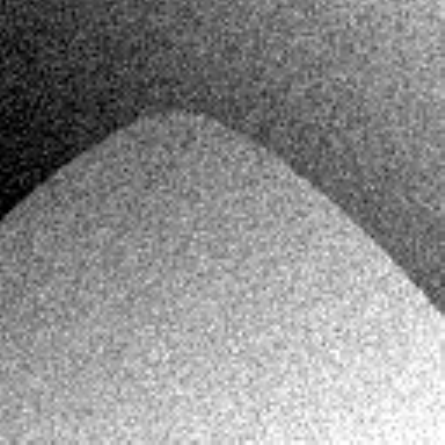}
\\ 
\includegraphics[width=5.8cm]{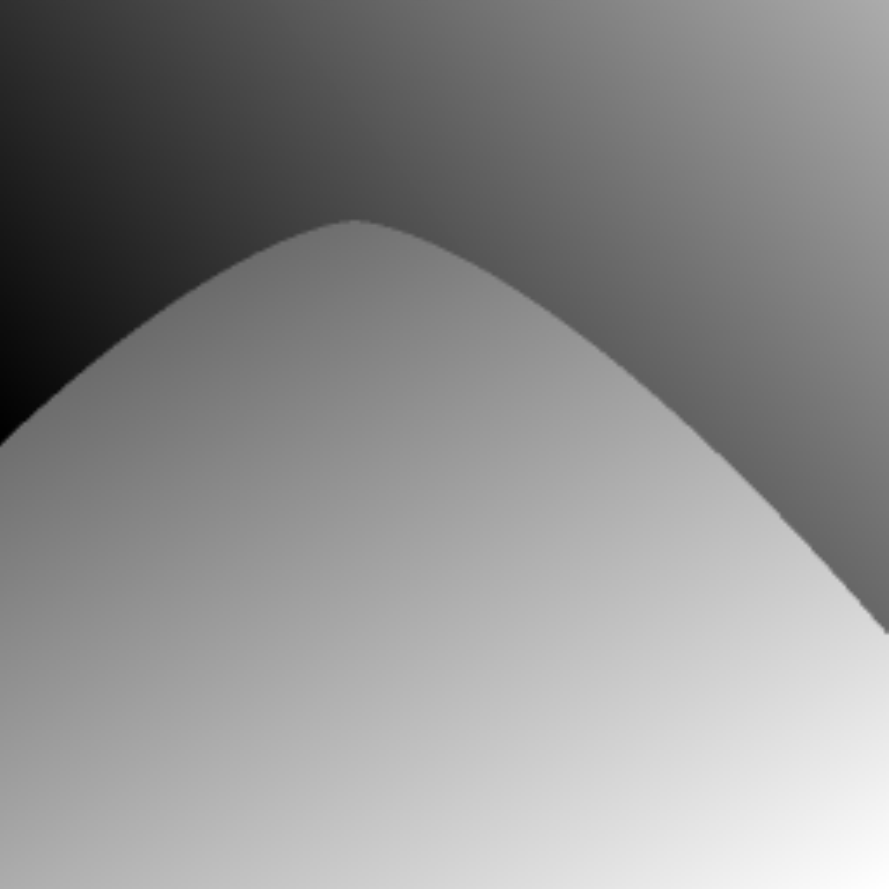}
\quad & 
\includegraphics[width=5.8cm]{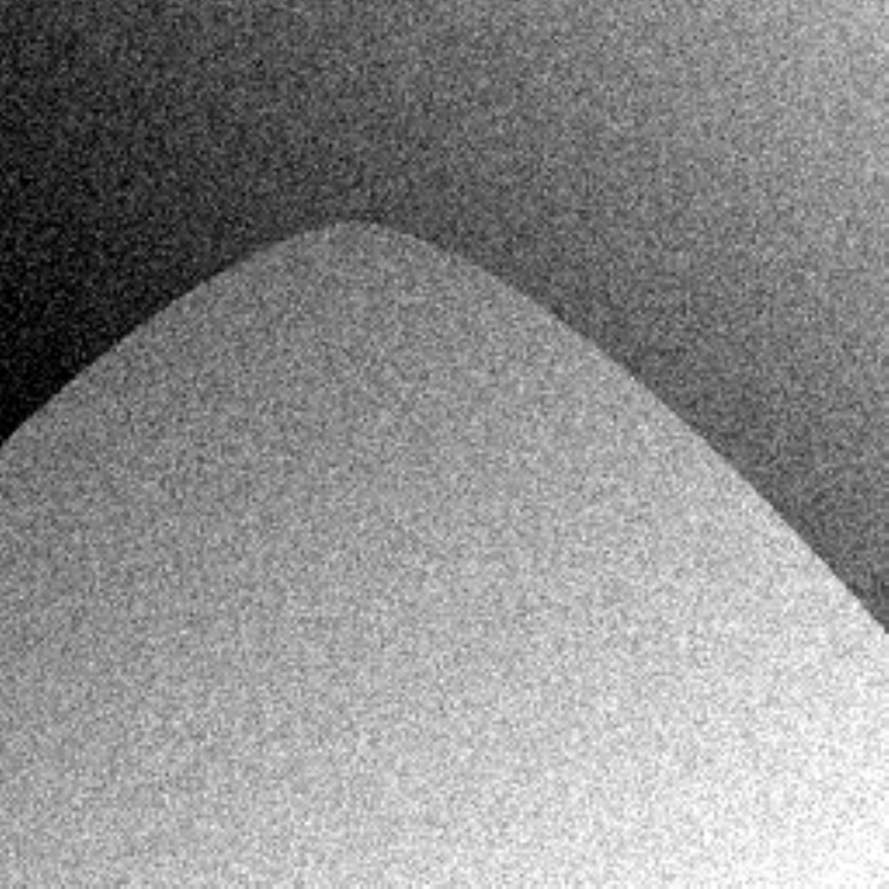}
\end{tabular}
\caption{Left: $\delta^n f$, for   $n=64,128,256$, respectively, where $f$ is a horizon function, according to Definition \ref{def:horizon_functions}. Here, the horizon boundary is in $\mathscr{C}^\alpha((0,1))$ and $\alpha = 1.5$. Right: Respective noisy images $\delta^n f+ \xi^n$.}
\label{fig:original_horizons}
\end{figure}

\subsubsection*{Wedgelets and approximations}
 
We first recall some approximation results for wedgelets.  They stem from  \fcite{Donoho:1999} and \fcite{WillettNowak:2003}. Since we are not working with the same discretisation we rewrite  them  for the continuous setting and provide elementary self-contained proofs. The discussion of the  discretisation is postponed to Section \ref{sec:wedgelets_consistency}. 
We start with the definition of horizon functions, like in \fcite{Donoho:1999}. 

\begin{definition}[Horizon functions]\label{def:horizon_functions}  Let $\alpha \in (1,2]$  and  $h \in
\mathscr{C}^{\alpha}([0,1]) $ if $\alpha <2$ or $\mathscr{C}^{1,1}([0,1])$ if $\alpha=2$. Let further $f$ be a bivariate function  which is piecewise constant relative 
to the partition of $[0,1]^2$ in an upper and a lower part induced by  $h$:
$$
   f(x,y) = \begin{cases}  c_1 \text{ if } y \leq h(x), \\ c_2 \text{ if } y > h(x), \end{cases}
$$
with real numbers  $c_1$ and $c_2$. Such a  function  is called an $\alpha$-\emph{horizon function}; the set of such functions will be 
denoted by $Hor^\alpha ([0,1]^2)$. $h$ is called the \emph{horizon boundary} of $f$.
\end{definition}

Discretisation at various levels of a typical horizon function is plotted 
Fig.~\ref{fig:original_horizons}, left column. In the right column the respective noisy versions are shown. 

\begin{lemma} \label{lem:Donoho} Let $\alpha \in [1,2]$ and  $f \in Hor^\alpha ([0,1]^2)$ 
with boundary function $h$. Then there are $C,C^\prime>0$ 
- independent of $k$ - and for each $k$ a continuous wedge partition 
$\mathscr{W}_{k}$ of the unit square $[0,1]^2$, such that
$|\mathscr{W}_{k}| \leq  C^\prime  k $ and 
$$
  \|f-f_k\|_{L^2([0,1]^2)} \leq  \frac{C}{k^{\alpha/2}},
$$
where $f_k$ is the $L^2$-projection of f on the space of piecewise constant functions relative to
the wedge partition $\mathscr{W}_k$.
\end{lemma}

\begin{proof} Let us first approximate the graph of $h$ by linear  pieces. 
We consider the uniform partition induced  by $x_i = i/k$.  We denote 
by $S_k(h)$ the continuous linear spline interpolating  $h$ relatively to  the uniform subdivision:
$$
   S_k(h)(x) = h(x_i) + (x-x_i) \left(\frac{h(x_{i+1}) - h(x_i)}{x_{i+1}-x_i} \right) \mbox{ for } i= 0,\hdots,k-1 \mbox{ and } x \in I_i
$$
where $I_i = [x_i,x_{i+1}]$. Therefore, we have
\begin{equation} \label{eq:1}
   \left|h(x)-S_k(h)(x) \right| =  \left| h(x) -h(x_i) - \frac{h(x_{i+1}) - h(x_i)}{x_{i+1}-x_i}    (x-x_i )
   \right| \; \text{ for each } x \in I_i.
\end{equation}
Since $h^\prime \in \mathscr{C}^{0,\alpha-1}([0,1])$, there exists  $C>0$  such that
$$
  \left| \frac{h(x_{i+1}) - h(x_i)}{x_{i+1}-x_i} - h^\prime(x_i) \right| \leq 
  C |x_{i+1}-x_i|^{\alpha-1} = \frac{C}{k^{\alpha-1}}.
$$
This implies that
$$
\left|h(x)-S_k(h)(x)\right| = \left| h(x) -h(x_i) - \left(h^\prime(x_i)
+O \left(\frac{1}{k^{\alpha-1}} \right) \right) (x-x_i)\right| \mbox{ for } x \in I_i.
$$
On the other hand,
$$
    h(x) = h(x_i) + h^\prime(x_i) (x-x_i) +O (|x-x_i|^{\alpha}).
$$
Hence, Equation \eqref{eq:1} can be
rewritten  as 
$$
    \left|h(x)-S_k(h)(x) \right| =  O(|x-x_i|^{\alpha }) +  O\left(\frac{1}{k^\alpha}\right)
$$
and there is a constant $C>0$ (independent of $k$) such that
\begin{displaymath} \label{eq:horizon_approx_estimate}
  \left\|h-S_k(h)\right\|_{L^\infty([0,1])} \leq \frac{C}{k^{\alpha}}.
\end{displaymath}
Now we will use this estimate
to derive error bounds
for the optimal wedge representation. As a piecewise approximation of $f$ we propose  
$$
   f_k(x,y) = \begin{cases}  c_1 \text{ if } y < S_k(h)(x); \\ c_2 \text{ if } y > S_k(h)(x). \end{cases}
$$
 We bound the error 
 by the area between the horizon $h$ and its
piecewise affine reconstruction:
\begin{align*} \label{eq:L_infty_estimate}
  \|f-f_k\|_{L^2([0,1]^2)} &\leq |c_1-c_2| \left( \int_0^1
  |h(x)-S_k(h)(x)|\,dx \right)^{1/2} \\ &\leq |c_1-c_2|
  \left(\|h-S_k(h)\|_{L^\infty([0,1])} \right)^{1/2}
  \leq \frac{C}{k^{\alpha/2}}.
\end{align*}



\noindent It remains to bound the size of the minimal continuous wedgelet partition $\mathscr{W}_k$, such that
$f_k \in \mathscr{F}_{\mathscr{W}_k}$. A proof  is given
in Lemma 4.3 in \fcite{Donoho:1999}; it uses $h \in
\mathscr{C}^1([0,1])$.
\end{proof}

\begin{bemerkung} For an arbitrary horizon function, the approximation
rates obtained by non-linear wavelet approximation (with sufficiently smooth wavelets)  can not be better  than  
$$
  \|f-f_k\|_{L^2([0,1]^2)} = O \left( \frac{1}{k^{1/2} }\right),
$$
where $f_k$ is the non-linear $k$-term wavelet approximation of $f$.
This means that for such a function the asymptotical behaviour in terms of approximation 
rates is strictly better for wedgelet decompositions than for wavelet decompositions. For a discussion on this topic, see Section 1.3 in \fcite{CandesDonoho:2002}.
\end{bemerkung}

Piecewise constant wedgelet representations are limited by 
the degree 0 of the regression polynomials on each wedge. This
is reflected by the choice of the horizon functions which are not
only piecewise smooth but even piecewise constant. A similar  phenomenon arises also
in the case of approximation by Haar wavelets. 

 \fcite{WillettNowak:2003} extended the regression model to piecewise linear functions on each leaf of the wedgelet
partition and called the according representations \emph{platelets}. 
This allows for an improved approximation rate for larger classes of piecewise smooth functions.   

Let $h$ be a function in $\mathscr{C}([0,1])$. We define the two subdomains $S^+$ and $S^-$, respectively, as the hypograph and the epigraph of $h$ restricted to $(0,1)^2$. In other words:
\begin{equation}
\label{eq:subdomains}
S^+= \left\{  (x,y) \in (0,1)^2  \; | \; y > h(x)\right\},  \;
S^-= \ag (x,y) \in (0,1)^2 \; | \; y < h(x) \ad. 
\end{equation}

Let us introduce the following generalised class of horizon
functions:
\begin{equation} \label{eq:ghorizon}
  Hor_1^{\alpha}([0,1]^2) := \{f:[0,1]^2 \rightarrow \R | \; f|_{S^+} \text{ and } f|_{S^-} \in \mathscr{C}^\alpha(S^\pm), h \in \mathscr{C}^\alpha([0,1])    \}.
\end{equation}
The following result from  \fcite{WillettNowak:2003} gives 
approximation rates by platelet approximations for $Hor^{\alpha}$. 


\begin{proposition} \label{prop:platelets_approx} Let $f \in
Hor_1^{\alpha}([0,1])$ for $1 < \alpha \leq 2$. 
Then the $k$-term platelet approximation error $h_k$
satisfies
\begin{equation} \label{eq:platelets_approx}
   \|f-h_k\|_{L^2([0,1]^2)}= O\left(\frac{1}{k^{\alpha/2}}\right).
\end{equation}
\end{proposition}

\begin{proof} A sketch of the proof  is given by the following two steps:
(1) the boundary between the two areas is approximated uniformly
like in \fcite{Donoho:1999}; (2) in the rest of the areas we use also
uniform approximation with dyadic cubes, together with the
corresponding H\"{o}lder bounds. The partition generated 
consists of squares  of sidelength at least $O(1/k^{1/2})$.
There are at most $O(k)$ such areas.
\end{proof}

\subsubsection*{Wedgelets and consistency}
\label{sec:wedgelets_consistency}

Now we apply the continuous approximation results to the consistency problem of the wedgelet estimator based on the above discretization. 
Note that, due to our specific discretization, the arguments below differ from those in  \fcite{Donoho:1999}.

Two ingredients are needed: pass over  to a suitable discretisation and bound the number of generated discrete wedgelet partitions polynomially in $n$, in order to apply the general consistency results. Let us first state a discrete approximation lemma:

\begin{lemma} \label{lem:wedgelets_discrete_approx} Let $f $ be an $\alpha$ horizon function in $Hor_1^\alpha$ with $1<\alpha<2$. 
There is $C>0$ such that for all $k \leq n \in \N$,  there is
 $(\mathscr{P}^n_k,h^n_k,) \in \mathfrak{S}^n$, such that 
$|\mathscr{P}^n_k| \leq k$ and which satisfies
\begin{equation}
\label{eq:wedgelet_approximation_discretisation}
\|f-\iota^n h^n_k\| \leq C \cdot \left( \frac{1}{k^{\alpha/2}} +    
 \frac{k^{1/2}}{n^{1/2}} \right).
\end{equation}
\end{lemma}

\begin{proof}  
The triangular inequality yields the following decomposition of the error
$$
 \|f-\iota^n \delta^n h_k\| \leq     \|f-    h_k \|+ \| h_k - \iota^n \delta_n h_k\|.
$$
The first term may be approximated by (\ref{eq:platelets_approx}), whereas the second term corresponds to the discretisation. Let us estimate the error induced by discretisation. 


One  just has to split $[0,1)^2$ into $J_n^k$, 
 the union of those squares in $\mathscr{Q}_n$ which do not intersect the approximating wedge lines and $K_n^k$ the union of such squares meeting the approximating wedge lines. We obtain the following estimates:
$$
   \|h_k - \iota^n \delta^n h_k \|^2_{L^2(Q)} \leq  \frac{C}{n^2}  \mbox{ for any } Q \in K_n^k, \mbox{ and for some constant } C > 0. 
$$
Since there are at most  $C^\prime   k n $ such squares, for some constant $C^\prime$ not depending on $k$ and $n$, this implies that
$$
   \|h_k - \iota^n \delta^n h_k \|^2_{L^2(K_n^k)} \leq \frac{C  kn}{n^2} = \frac{C }{n}
\text{        and       }
 \|h_k - \iota^n \delta^n h_k \|^2_{L^2(J_k^n)} \leq \frac{C k}{n},
$$
where $C>0$ is a constant.
Taking $h_k^n = \delta^n h_k$ completes the proof.
\end{proof}



Finally, the following lemma provides an estimate of the number of fragments in $\mathscr{R}^n$.

\begin{lemma}  \label{lemma:counting_wedgelets} There is a constant $C>0$ such that for all $n \in \M$ the number
$|\mathscr{R}^n|$ of fragments used to form the wedgelet partitions is bounded as follows:
$$
   |\mathscr{R}^n| \leq C n^4.
$$ 
\end{lemma}

\begin{proof}  
In a dyadic square of size $j$, there are at most $j^4$ possible digital lines. 
For dyadic $n \in \M$  one can write  $n=2^J$ and therefore we have
$$
     |\mathscr{R}^n| \leq \sum_{i=0}^J 2^{2(J-i)} \cdot 2^{2 \cdot 2i}  = n^2 \sum_{i=0}^J 2^{2i} = n^2 \cdot \frac{2^{2J+2}-1}{2^2-1} \leq C \cdot n^4 \mbox{ for some constant } C > 0.
$$ 
This completes the proof.
\end{proof}

Note that the discretisation of the continuous approximation $h_k$ leads to a wedgelet partition composed of fragments in $\mathscr{R}^n$. Therefore, combination of the Lemmata \ref{lemma:counting_wedgelets} and \ref{lem:wedgelets_discrete_approx} yields:

\begin{theorem} \label{th:2d_consistency_wedgelets} On $S^n = \{1,\,\ldots \, , n\}^2$ the following holds: Let $\alpha$ such that $1<\alpha<2$
and $f$ be an $\alpha$ horizon function in $Hor_1^\alpha([0,1]^2)$ with $1<\alpha<2$ and suppose that $\gamma_n$ satisfy (H\ref{hyp rate assumptions}.3).
 Assume further that the noise variables $\xi_s^n$ from Section \ref{suse The Regression Model} satisfy \eqref{eq MainTailEstimate}. 
Then 
\begin{equation} \label{eq:2D_consistency_wedgelets}
   \| \hat{f}^n_{\gamma_n} -f\|^2 = O \left( \gamma_n^{\frac{\alpha}{\alpha+1}}\right) + 
   O\left(n^{-\frac{ \alpha}{\alpha+1}}\right), 
\hbox{  for almost all }\omega \in\Omegait,
\end{equation}
where $\hat{f}^n_{\gamma_n} $ is the wedgelet-platelet estimator.
%
\end{theorem}

\begin{bemerkung} \label{rem:wedgelets_discretisation}Choosing $\gamma_n$ of the order $\ln n /n^2$, estimate $\eqref{eq:2D_consistency_wedgelets}$ reads 
\begin{equation} \label{eq:2D_consistency_wedgelets2}
\| \hat{f}^n_{\gamma_n} -f\|^2 = O \left( \frac{(\ln n)^{\frac{2 \alpha}{\alpha+1}} }{n^{\frac{2\alpha}{\alpha + 1}}} \right) +  
O\left(\frac{1}{n^{\frac{ \alpha}{\alpha+1}}}\right) 
\hbox{  for almost all }\omega \in\Omegait.
\end{equation}
Whereas the left term on the right-hand size consists of the best compromise between approximation and noise removal, the right term on the right-hand size corresponds to the discretisation error. Note that, in contrast to the $1D$-case the discretisation error dominates, as soon as $\alpha > 1$. This is related to the piecewise constant nature of our discretisation operators. In concrete applications, this may prove to be a severe limitation to the actual quality of the estimation. Up to this discretisation problem, 
the decay rates given by \eqref{eq:2D_consistency_wedgelets2} are the usual optimal rates for the function class under consideration.
\end{bemerkung}

\begin{figure}[t!!]
\begin{tabular}{cc}
\includegraphics[width=5.8cm]{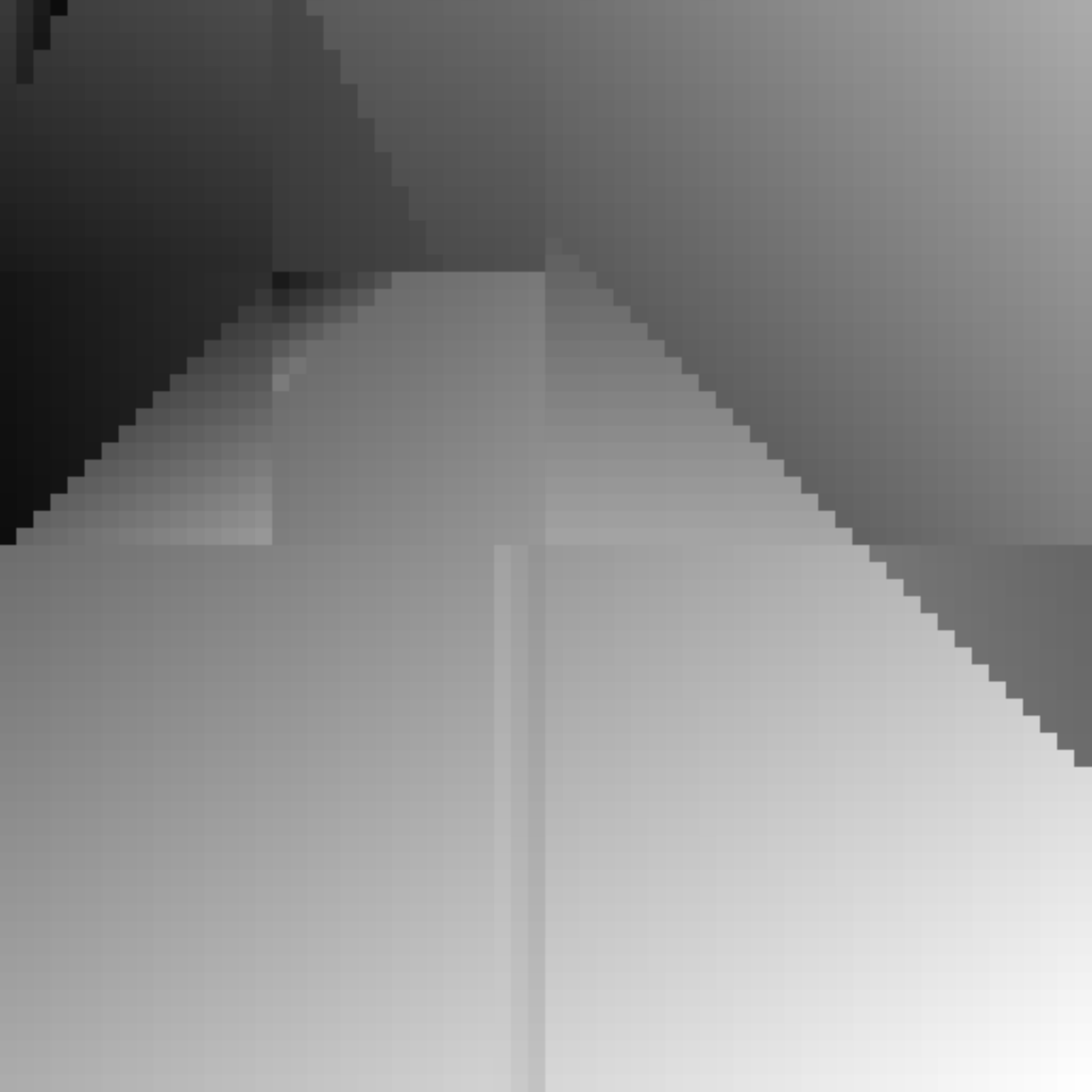}
\quad & 
\includegraphics[width=5.8cm]{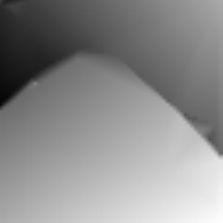}
\\
\includegraphics[width=5.8cm]{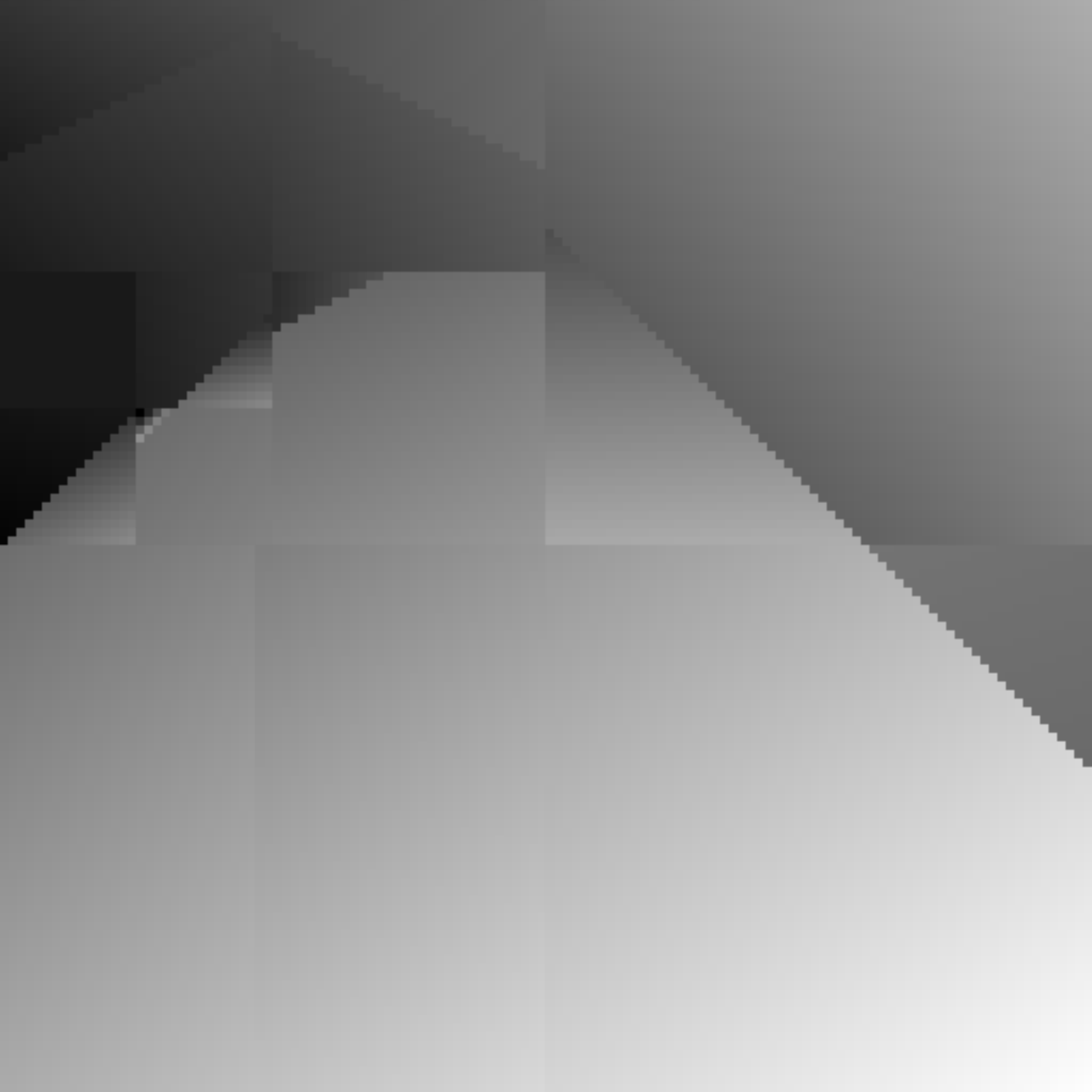}
\quad &
\includegraphics[width=5.8cm]{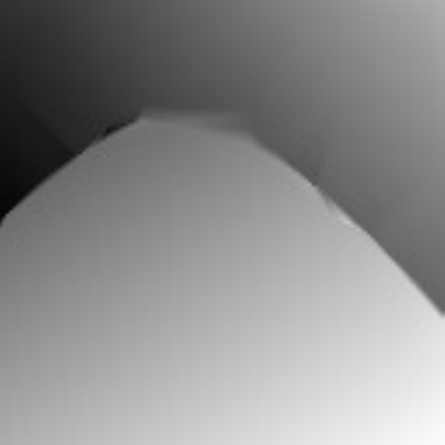}
\\ 
\includegraphics[width=5.8cm]{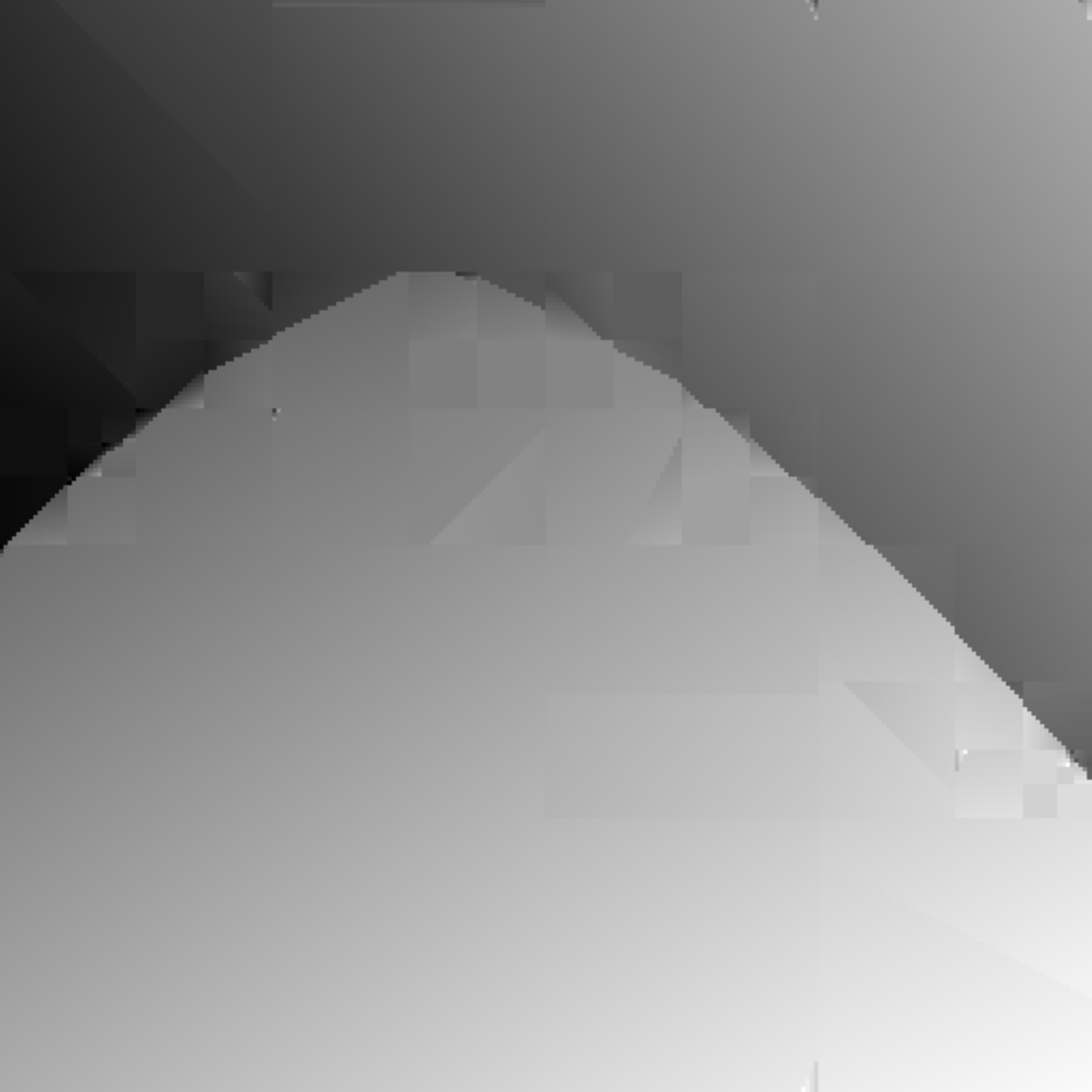}
\quad & 
\includegraphics[width=5.8cm]{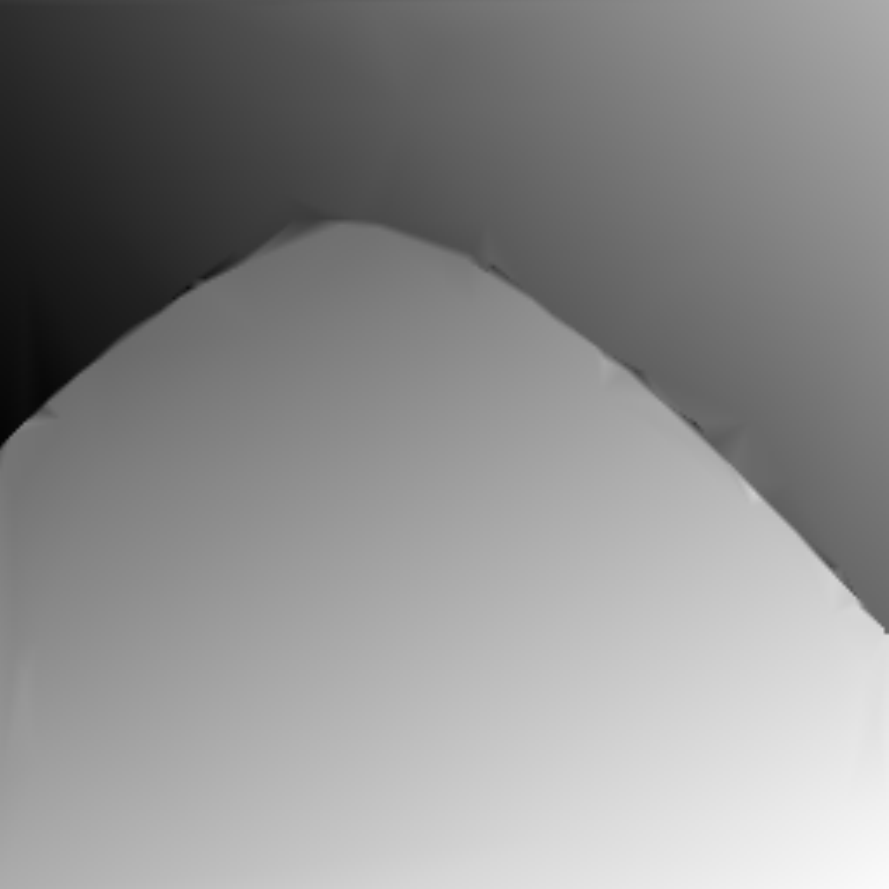}
\end{tabular}
\caption{Estimators of the noisy images of Fig~\ref{fig:original_horizons}. Left: piecewise linear wedgelet estimator. Right: piecewise linear and continuous Delaunay estimators.}
\label{fig:estimators}
\end{figure}

On the left column of Fig.~\ref{fig:estimators}, wedgelet estimators for a typical noisy horizon function are shown.  

\subsection{Triangulations}

Adaptive triangulations have been used since the emergence of early finite element methods to approximate solutions of elliptic differential equations. They have been also used in the context of image approximation; we refer to \fcite{DemaretIske:2010} for an account on recent triangulation methods applied to image approximation. The idea to use discrete triangulations leading to partitions based on a polynomially growing number of 
triangles has been proposed in \fcite{Candes:2005} in the context of piecewise constant functions over triangulations.  In the present example, we deal with a different approximation scheme, where the triangulations are Delaunay triangulations and were the approximating functions are continuous linear splines.
One key ingredient is the use of  recent approximation results, \fcite{DemaretIske:2011},  that show the asymptotical optimality of  approximations based on Delaunay triangulations having at most $n$ vertices.  Due to this specific approximation context, a central ingredient  for the proof of the consistency is a suitable  discretization scheme, which still preserves the approximation property. 


\subsubsection*{Continuous and discrete triangulations}
\label{sec:continuous_discrete_triangulations}
Let us start with some definitions. We begin with triangulations in the continuous settings:
\begin{definition} 
\label{def:triangulation} 
A \emph{conforming triangulation} $\mathcal{T}$ of the domain $[0,1]^2$ is a finite set $\{T\}_{T \in \mathcal{T}}$
of closed triangles $T \subset [0,1]^2$ satisfying the following conditions.
\begin{itemize}
\item[(i)] The union of the triangles in $\mathcal{T}$ covers the domain $[0,1]^2$;
\item[(ii)]
for each pair $T,T^\prime \in \mathcal{T}$ of distinct triangles,
the intersection of their interior is empty;
\item[(iii)]  any pair of two distinct triangles in $\mathcal{T}$ intersects at most in 
one common vertex or along one common edge.
\end{itemize}
We denote the set of (conforming) triangulations by $\mathscr{T}([0,1]^2)$. We will use the term
\emph{triangulations} for conforming triangulations.
\end{definition}
Accordingly we define the following discrete sets, relatively to
 partitions  $\mathscr{Q}_k=\{[(i-1)/k,i/k)  \times [(j-1)/k,j/k): i,j=1,\cdots,k\}$ of $[0,1)^2$ into $k$ squares each of side length $1/k$.

For $a,b \in [0,1]^2$ we denote by $[a,b]$ the line segment with endpoints $a$ and $b$.

\begin{definition} \label{def:discrete_triangulation}
For a triangle $T \subset [0,1]^2$, with vertices $a,b$ and $c$, we define  the following discrete sets:
\begin{itemize}
\item[(i)] for each $p\in \{a,b,c\}$ the square $Q \in \mathscr{Q}_n$ such that $Q \ni p$  is called a {\em discrete vertex} of $T$;
\item[(ii)] for each edge $e \in \{[ab],[bc],[ca]\}$, the set of squares $Q \in \mathscr{Q}_n$ such that $Q \cap e \neq \emptyset$ and $Q$ is not a discrete 
vertex is called a {\em discrete (open) edge} of the triangle $T$;
\item[(iii)] the set of squares $Q \in \mathscr{Q}_n$ such that $Q \cap T \neq \emptyset$ and $Q$ is neither a discrete vertex nor belongs to a  discrete open edge is called a {\em  discrete open triangle}.
\end{itemize}
\end{definition}


\subsubsection*{Piecewise polynomials functions on triangulations}
\label{sec:Piecewise_polynomials_triangulations}

We take $S^n = \{1,\cdots,n\}^2$ and the set of fragments $\mathscr{R}^n$ is given as the set of discrete vertices, open edges  and open triangles 
$$\mathscr{R}^n =S^n\cup \left\{ ([ab]) : a,b \in S^n \right\} \cup \left\{ ([abc]) : a,b,c \in S^n \right\}.$$ 
 We let $\mathfrak{P}^n$ then be  the collection of partitions of  $S^n$ into discrete triangles, obtained from a continuous triangulations, and assuming that there is a rule deciding to which triangle discrete open segments and discrete vertices belong. Each such discrete triangle is then  the  union of elementary  digital sets in $\mathscr{R}^n$. We remark that  $|\mathscr{R}^n|=n+n(n-1)/2 + n(n-1)(n-2)/6$ and therefore $|\mathscr{R}^n| \sim n^3/6$.
 Like in the one-dimensional case, as described in Section~\ref{sec:1D_interval_partitions},  we choose the finite dimensional linear subspace $\mathscr{F}_p \subset L^2([0,1))$ of polynomials of maximal degree $p$. The induced segmentation  classes
$\mathfrak{S}^n(\mathfrak{P}^n,\mathfrak{F}^n)$ consist of piecewise polynomial functions relative to partitions in $\mathfrak{P}^n$.  

We have the following approximation lemma
\begin{lemma} \label{lem:regular_triangulations_approx} Let $f \in \mathscr{C}^\alpha([0,1]^2)$, with $p<\alpha<p+1$. There is $C>0$
such that for all $k \leq n \in \N$, 
there is
 $(\mathscr{P}^n_k,h^n_k) \in \mathfrak{S}^n$, such that 
$|\mathscr{P}^n_k| \leq k$ and which satisfies
\begin{equation}
\label{eq:approximation_discretisation 2}
\|f-\iota^n h^n_k\| \leq C \cdot \left( \frac{1}{k^{\alpha/2}} +    
 \left( \frac{k}{n}\right)^{1/2} \right).
\end{equation}
\end{lemma}

\begin{proof} We first use classical aproximation theory which tells us the existence of a function $h_k:[0,1]^2 \mapsto \R $,  piecewise polynomial relatively to  a triangulation with $k$ triangles and such that the error on the whole domain is bounded by
$$
    \|f-h_k\| \leq \frac{C}{k^{\alpha/2}}.
$$ 
As in the 1-D case we split $[0,1)^2$ into the union $J_n^k$ of those squares in $\mathscr{Q}_n$ which do not meet the continuous triangulation, and 
  $K_n^k$ the set of such squares meeting the triangulation, i.e. which intersects with some edge of the triangulation.
For each small square $Q \in \mathscr{Q}_n$ and $Q \subset K_n^k$, the following estimate holds:
$$
   \|h_k - \iota^n \delta^n h_k\|^2_{L^2(Q)} \leq \frac{C}{n^2} \mbox{ for any } Q \in K_n^k, \mbox{ and some constant } C>0
$$
and there are at most $3 \cdot \sqrt{2} kn$ such squares. 
Altogether we obtain:
$$
   \|h_k - \iota^n \delta^n h_k\|_{L^2(K_n^k)} \leq \frac{C k^{1/2}}{n^{1/2}}, \mbox{ for some constant } C > 0.
$$
Now for each square $Q \in \mathscr{Q}_n$ and $Q \subset J_n^k$,  an argumentation similar to that  in  the $1D$-proof yields 
$$
   \|h_k -\iota^n \delta^n h_k\|_{L^2(J_k^n)} \leq \frac{C}{n}.   
$$
This completes the proof.\end{proof}

Due to  Lemma \ref{lem:regular_triangulations_approx},
(\ref{eq Bedingung mit F}) is satisfied:   a function in $\mathscr{C}^\alpha$ satisfies (\ref{eq Bedingung mit F}) with $\rho = 1/2$ and $F_n = n^{1/2}$
and therefore 
Theorem \ref{th convergence rates} applies. 

\subsubsection*{Continuous linear splines}
\label{sec:Piecewise_polynomials_continuous_delaunay_triangulations}

We turn now to the more subtle  case of continuous linear splines on Delaunay triangulations.
Anisotropic Delaunay triangulations have been recently applied successfully to  the design of a full image compression/decompression scheme,
\fcite{DemaretDynIske:2006}, \fcite{DemaretIskeKhachabi:2009}.  We  apply such triangulation schemes in the context of image estimation.  


We first introduce the associated function space in the continuous setting. We restrict the discussion to the case of piecewise affine functions, i.e. $p=1$.

\begin{definition} \label{def:continuous_linear_splines_triangulations} 
Let $\mathcal{T} $ be a conforming triangulation of $[0,1]^2$. Let
$$
  \mathcal{S}^0_{\mathcal{T}} = \left\{    f \in \mathscr{C}\left([0,1]^2\right) \; {\bf :} \; f \big|_T \in \mathscr{F}_1,
  T \in \mathcal{T} \right\},
$$
be the set of piecewise affine and continuous functions on $\mathcal{T}$.
\end{definition}

The following  piecewise smooth functions generalise the horizon functions from (\ref{eq:ghorizon}).

\begin{definition}
\label{def:generalised_horizon}  Let $\alpha \in (1,2)$ and $g \in \mathscr{C}^\alpha([0,1])$. Let $S^+$ and $S^-$ be two subdomains defined as in (\ref{eq:subdomains}). A \emph{generalised $\alpha$-horizon function} is an element of the set 
$$
     \mathscr{H}^{\alpha,2}([0,1]^2) := \ag f \in L^2([0,1]^2)   \; | \;  f|_{S^+},f|_{S^-} \in  W^{\alpha,2}(S^{\pm}) \ad 
$$
where $W^{\alpha,2}(S^\pm)$ is the Sobolev space of regularity $\alpha$ relative to the $L^2$-norm on $S^\pm$.
 \end{definition}

In order to obtain convergence rates of the triangulation-based estimators for this class of functions we  need the following recent result,  Thm.4 in \fcite{DemaretIske:2011}:
\begin{theorem} \label{th:Delaunay_approx}  Let $f$ be an $\alpha$-horizon function in $Hor_1^\alpha$, with $\alpha \in (1,2)$, such that $ f|_{S^\pm} \in W^{\alpha,2}(S^\pm)$.  Then there is $C>0$, such that for all $k \in \N$ there is a Delaunay triangulation $\mathcal{D}_k$ with
$$
    \|f - \pi_{\mathcal{S}^0_{\mathcal{D}_k}} f\|_{L^2([0,1]^2)} \leq \frac{C}{k^{ \alpha/2 }}.
$$
\end{theorem}

Using arguments as in the proof of Lemma \ref{lem:regular_triangulations_approx}, we  obtain the following lemma:
\begin{lemma} \label{lem:discrete_Delaunay_triangulations_approx} Let $f \in \mathscr{H}^{\alpha,2}([0,1]^2)$, with $1<\alpha<2$ there is $C>0$
such that for all $k \leq n \in \N$, 
there is
 $(\mathscr{P}^n_k,h^n_k)$, such that 
$\mathscr{P}^n_k \in \mathfrak{P}^n$ is a discretisation of a continuous Delaunay triangulation $\mathcal{D}_k$, $|\mathscr{P}^n_k| \leq k$,  $h_k^n = \delta^n h_k$, where $h_k \in \mathcal{S}^0_{\mathcal{D}_k}$ and which satisfies
\begin{equation*}
\|f-\iota^n h^n_k\| \leq C \cdot \left( \frac{1}{k^{\alpha/2}} +    
\frac{k^{1/2}}{n^{1/2}} \right).
\end{equation*}
\end{lemma}

The previous machinery  cannot be applied directly without an explanation: since we are dealing with the space of continuous linear splines, our scheme is not properly a projective $\mathscr{F}$-segmentation class.  However, for each fixed partition, $\mathscr{P} \in \mathfrak{P}$ with elements  in $\mathscr{R}^n$, $\mathcal{S}^0_{\mathcal{T}}$ a subspace of $\mathscr{F}_{\mathscr{P}}$. Observe that all arguments in Lemma \ref{lem:fundamental_inequality} remain valid if we replace $\mathscr{F}_{\mathscr{P}}$ by subspaces and consider also the minimisation of the functional $H_\gamma^n$ over functions in these subspaces. We can therefore apply Theorem \ref{th convergence rates} to obtain the equivalent of Theorem \ref{th:1d_consistency}.

\begin{theorem} \label{th:2d_consistency_triangulation} Let $1<\alpha<2$
and let $f$ be a generalised horizon function in $\mathscr{H}^\alpha([0,1]^2)$. Let further assume that
noise in \eqref{eq regression model}  satisfies \eqref{eq MainTailEstimate} and that
 $\gamma_n$ satisfy (H\ref{hyp rate assumptions}.3).
Then 
\begin{equation} \label{eq:2D_consistency_triangulation}
   \| \hat{f}^n_{\gamma_n} -f\|^2 = O \left( \gamma_n^{\frac{\alpha}{\alpha+1}}\right) + O\left( n^{-\frac{ \alpha}{\alpha+1}}\right) 
\hbox{  for almost all }\omega \in\Omegait,
\end{equation}
where $\hat{f}^n_{\gamma_n}$ is the Delaunay estimator. 
%
\end{theorem}

\begin{proof} We check  the assumptions in Theorem \ref{th convergence rates}. 
Since $|\mathscr{R}^n| $ is of the order $ (n^2)^3$,
 Hypothesis  {\em (H\ref{hyp rate assumptions}.1)} holds with $\kappa = 3$. Hypothesis {\em (H\ref{hyp rate assumptions}.2)}
and  {\em (H\ref{hyp rate assumptions}.3)} were required separately. Finally, (\ref{eq Bedingung mit F}) holds with $\varrho=1/2$ and $F_n=n^{1/2}$
by Lemma  \ref{lem:discrete_Delaunay_triangulations_approx}. This completes the proof.
\end{proof}

\begin{bemerkung} \label{rem:Triangulation_discretisation} Similarly to Remark~\ref{rem:wedgelets_discretisation} and choosing $\gamma_n$ of the order $\ln n /n^2$, estimate $\eqref{eq:2D_consistency_triangulation}$ reads 
$$
\| \hat{f}^n_{\gamma_n} -f\|^2 = O \left( \frac{(\ln n)^{\frac{2 \alpha}{\alpha+1}} }{n^{\frac{2\alpha}{\alpha + 1}}} \right) +  O\left(\frac{1}{n^{\frac{ \alpha}{\alpha+1}}}\right) 
\hbox{  for almost all }\omega \in\Omegait.$$
The discussion  of Remark~\ref{rem:wedgelets_discretisation}  can be easily adapted to the case of estimation by triangulations. 
\end{bemerkung}

On the right column of Fig.~\ref{fig:estimators},  estimators by Delaunay triangulation are shown, for the same noisy horizon function as in the wedgelet case. 

The rates in Theorem \ref{th:2d_consistency_triangulation} are, up to a logarithmic factor, similar to the minimax rates obtained in \fcite{CandesDonoho:2002} with curvelets for $\alpha = 2$ and more recently  in  \fcite{DossalMallatLePennec:2011} with bandelets for general $\alpha$. This is in contrast to isotropic approximation methods, e.g. shrinkage of tensor product wavelet coefficients, which only attain the rate for $\alpha = 1$.

\section{Appendix}\label{se appendix}
We are going now to supply the proof of Theorem \ref{hypsubgaussian}.  

\begin{proof}[ Proof of Theorem \ref{hypsubgaussian}] Suppose that \emph{(b)} holds.
Theorem 1.5 in \cite{BuldyginKozachenko:2000} gives
\begin{displaymath}
\mathbb{P}\left(\left|\sum_{s\in S_n}\mu_s\xi_s^n\right|\ge c\right)\le2\cdot\exp\left(-\displaystyle\frac{c^2}{2\sum_{s\in S_n}\tau^2(\mu_s\xi_s^n) }\right).
\end{displaymath}
$\tau$ is a norm and therefore  $\tau^2(\mu_s\xi_s^n)= \mu_s^2\tau^2(\xi_s^n)$. Because of (\ref{eq defBeta}), the inequality (\ref{eq MainTailEstimate}) holds and hence\emph{(a)}. 
Part \emph{(b)} follows from the following two Lemmata \ref{la SGVolkmar1} and  \ref{la SGVolkmar2}. In fact,  for $s\in S^n$ and $\mu_{s'}=\delta_{s,s'}$ the inequality (\ref{eq MainTailEstimate}) boils down to
$
\mathbb{P}(|\xi_s|\ge c)\le 2\exp(-c^2/\beta)
$
and the lemmata apply.
\end{proof}

The missing lemmata read:

\begin{lemma}\label{la SGVolkmar1}
Let $\xi$ be a random variable with
$
\mathbb{P}(|\xi|\ge c)\le C\cdot\exp(-c^2/\beta)
$, $\beta>0$.
Then
\begin{displaymath}
\mathbb{E}(\exp(t\xi^2))\le 1+ Ct/(\beta^{-1}-t)\hbox{   whenever  } |t| < 1/\beta.
\end{displaymath}
\end{lemma}
\begin{proof}
Let $\varrho$ be the distribution of $|\xi|$. With $b=1/\beta$ one computes
\begin{eqnarray*}
&\mathbb{E}(e^{t\xi^2})-1=\int_0^\infty e^{t x^2}\,d\varrho(x)-1=\int_0^\infty\int_0^x 2 t y e^{ty^2}\,dy\,d\varrho(x)=\int_0^\infty 2 t y e^{ty^2}\int_y^\infty \,d\varrho(x)\,dy\\
&=\int_0^\infty 2 t  y e^{ty^2}\mathbb{P}(|\xi|\ge y)\,dy\le 2 C t  \int_0^\infty  y e^{(t-b) y^2}\,dy=C t/(b-t)\hbox{     if }|t| < b.
\end{eqnarray*}
The proof is complete.
\end{proof}
\begin{lemma}\label{la SGVolkmar2}
Let $\alpha\ge 0$, $ \delta\ge 1$. Then there is $\beta' \in \mathbb{R}_+\cup\{\infty\} $ such that for all centred random variables $\xi$ with $\mathbb{E}(\exp(\alpha\xi^2))\le\delta$ the estimate
$\mathbb{E}(\exp (t\xi) )\le \exp( t^2/\beta')$ holds for every $t\in\mathbb{R}$.
\end{lemma}
A converse holds as well.
\begin{proof}Assume without loss of generality that $\alpha=1$.
Let us first  consider the case $|t|\ge2\ln^{1/2} \delta$.  Since 
$(\xi-t/2)^2\ge0$ one has $\exp(t\xi)\le \exp(t^2/4)\exp(\xi^2)$. Take expectations on both sides and use the assumption to get $\mathbb{E}(\exp(t\xi))\le \delta\exp(t^2/4)$. This implies
\begin{displaymath}
\mathbb{E}(\exp(t\xi))\le \exp(t^2/2)\hbox{   whenever }|t|\ge2\sqrt{\ln \delta}.
\end{displaymath}
Note that this estimate does not depend on the special variable $\xi$.

\noindent Let now 
$\label{equ t bounds}
|t|\le 2(\ln\delta)^{1/2}. 
$
The function $\varphi(t)=\ln\mathbb{E}(\exp(t\xi))$ is convex and hence $\varphi''(t)\ge 0$; furthermore 
\begin{equation}\label{equ conditions on phi}
\varphi(0)=0\hbox{  and  } \varphi'(0)=\mathbb{E}(\xi)=0.
\end{equation}
By the mean value theorem there is some $\vartheta(t)\in[0,1]$ such that
\begin{equation}
\varphi(t)=\varphi(0)+t\varphi'(0) + (t^2/2)\varphi''(\vartheta(t)t)\le (t^2/2)\max\{\varphi''(t):|t|\le 2\sqrt{\ln \delta}\}.
\end{equation}
Hence $1/\max(\max\{\varphi''(t):|t|\le2\ln^{1/2} \delta\},1)$ is a suitable scale factor for the  $\xi$ in question. 

\noindent We must finally remove           the dependency on moments of $\xi$ in
\begin{equation}\label{equ phi second derivative}
 \varphi''(t)=\left(\mathbb{E}(\xi^2\exp(t\xi))\mathbb{E}(\exp(t\xi))-\mathbb{E}^2(\xi\exp(t\xi))\right)  \left/\mathbb{E}^2(\exp(t\xi))\right..
\end{equation}
To this end let $\eta$ be an independent copy of $\xi$. Then the denominator 
becomes
\begin{displaymath}
D=\mathbb{E}(\xi^2\exp(t(\xi+\eta))-\mathbb{E}(\xi\eta\exp(t(\xi+\eta)))=\mathbb{E}((\xi-\eta)^2\exp(t(\xi+\eta))).
\end{displaymath}
With $(a-b)^2\le(a-b)^2+(a+b)^2=2(a^2+b^2)$ we arrive at
\begin{displaymath}
D\le 2\mathbb{E}(\xi^2\exp(t\xi))\mathbb{E}(\exp(t\xi)).
\end{displaymath}
By convexity of $\varphi$ and (\ref{equ conditions on phi}) one has $\varphi\ge 0$ and thus $\mathbb{E}(\exp(t\xi))\ge 1$. Furthermore, $\xi^4\le 2\exp(\xi^2)$. In view of the restriction on  $t$, the Cauchy-Schwartz inequality gives
\begin{displaymath}
\mathbb{E}^2\left(\xi^2\exp(t\xi)\right)
\le\mathbb{E}(\xi^4)\mathbb{E}(\exp(2t\xi))
\le2\cdot\mathbb{E}^2(\xi^2)\exp(t^2)
\le2\cdot\delta^2\cdot\delta^4
=2\delta^6.
\end{displaymath}
By Jensen's inequality  $\mathbb{E}(\exp(t\xi))\ge \exp(t\mathbb{E}(\xi))=\exp(t\cdot0)=1$. Hence
\begin{displaymath}
D\le 2^{3/2}\delta^3\mathbb{E}(\exp(t\xi))\le2^{3/2}\delta^3\mathbb{E}^2(\exp(t\xi)).
\end{displaymath}
Canceling out the numerator in (\ref{equ phi second derivative}) yields $\max\{\varphi''(t):|t|\le2(\ln(\delta))^{1/2}\}\le2^{3/2}\delta^3$ which completes the proof.
\end{proof}

\bibliographystyle{plainnat}
\bibliography{Consistencies}
\end{document}